\documentclass[a4paper,12pt,leqno]{article}
\usepackage{amsmath}
\usepackage{amsthm}
\usepackage{amssymb}
\usepackage{amscd}
\usepackage{graphicx, color}
\usepackage[all]{xy}

\newtheorem{thm}{Theorem}[section]
\newtheorem{definition}[thm]{Definition}
\newtheorem{lemma}[thm]{Lemma}
\newtheorem{crl}[thm]{Corollary}
\newtheorem{example}[thm]{Example}
\newtheorem{remark}[thm]{Remark}
\numberwithin{equation}{section}

\newcommand{\N}{\mathbb{N}}
\newcommand{\Z}{\mathbb{Z}}
\newcommand{\C}{\mathbb{C}}
\newcommand{\R}{\mathbb{R}}
\newcommand{\K}{\mathbb{K}}
\newcommand{\F}{\mathbb{F}}

\renewcommand{\P}{{\rm P}}

\newcommand{\RP}{\mathbb{R}{\rm P}}

\newcommand{\Map}{\mbox{{\rm Map}}}

\newcommand{\CP}{\mathbb{C}{\rm P}}

\newcommand{\dis}{\displaystyle}

\newcommand{\I}{\mbox{{\rm (i)}}}
\newcommand{\II}{\mbox{{\rm (ii)}}}
\newcommand{\III}{\mbox{{\rm (iii)}}}
\newcommand{\IV}{\mbox{{\rm (iv)}}}

\newcommand{\Po}{\mbox{{\rm Poly}}}
\newcommand{\po}{\mbox{{\rm Poly}}}
\newcommand{\pomod}[1]{\mathcal{P}oly^{#1,3}_{1}}

\title{\bf
Homotopy stability for spaces of triples of real polynomials
without common roots
}
\author{\bf Andrzej Kozlowski\footnote{%
Institute of Applied Mathematics and Mechanics,
University of Warsaw, Banacha 2, 02-097 Warsaw, Poland
(E-mail: akoz@mimuw.edu.pl)
}
\  and \ 
Kohhei Yamaguchi\footnote{%
Department of Mathematics,
University of Electro-Communications,  Chofu, Tokyo 182-8585, Japan
(E-mail: kohhei@im.uec.ac.jp)
}
}
\date{}
\begin{document}
\maketitle

\begin{abstract}
For each pair $(m,n)\not= (1,1)$ of positive integers
and a field $\F$ with algebraic closure $\overline{\F}$,
let
$\Po^{d,m}_n(\F)$ denote the space of
$m$-tuples
$(f_1(z),\cdots ,f_m(z))\in \F [z]^m$ of monic polynomials with coefficients in $\F$ of the same degree $d$ with no common roots in $\overline{\F}$ of multiplicity
$\geq n$.

These spaces were first explicitly defined and studied in an algebraic setting by
B. Farb and J. Wolfson in order to prove algebraic analogues of certain topological results of Arnold, Segal, Vassiliev, and others. 
They possess stability properties that have attracted considerable interest. 
For $\F=\C$ and $mn\geq 3$, our earlier work established homotopy
stability and determined the stable homotopy types of these spaces.
For $\F=\R$ we obtained analogous results when $mn\geq 4$ and
proved homological stability when $mn=3$.
In this paper we settle the case $(m,n)=(3,1)$: the space of triples of
monic real polynomials of the same degree having no common root.  We prove
homotopy stability and determine the stability range supplied by
our method.  
We also show that, in degrees $2$ and $3$, the distinction between homotopy
equivalence up to and through the stability dimension is essential.
We have also proved homotopy stability in the Arnold--Vassiliev case
$(m,n)=(1,3)$, for single monic real polynomials with no root of
multiplicity at least three in $\C$. The proof will appear in a separate paper.
\end{abstract}

\section{Introduction}\label{section: introduction}

\paragraph{1.1  Historical survey.}
The motivation of this paper comes from the work of 
B. Farb and J. Wolfson \cite{FW}. 
Inspired by the classical theory of resultants, they defined an algebraic variety $\po^{d,m}_n(\F)$.
In particular, they
computed various algebraic and geometric invariants of these varieties in order to address a conjecture
when $\F=\F_q$ (a finite field).
Moreover, for the case $\F =\C$,
the homotopy type of $\po^{d,m}_n(\F)$ has been extensively studied
by several mathematicians 
(e.g. \cite{Ar},  \cite{CCMM},   \cite{GKY2}, \cite{KY8}, \cite{KY9}, \cite{KY12}, \cite{Se},  \cite{Va}).
In this paper we investigate the space $\po^{d,m}_n(\F)$ in the case $\F=\R$.
For this purpose, recall the definition of the algebraic variety $\po^{d,m}_n(\F)$:
\begin{definition}
{\rm 
For each pair $(m,n)\not= (1,1)$ of positive integers and
a field $\F$ with its algebraic closure $\bar{\F}$, 
let $ \Po^{d,m}_n(\F)$ be the space of $m$-tuples
$(f_1(z),\cdots ,f_m(z))\in \F [z]^m$ 
of monic polynomials with coefficients in $\F$ of the same degree $d$ with no common root 
in $\overline{\F}$ of multiplicity $\geq n$. 
\par
Note that there is a homeomorphism
\begin{equation}\label{eq: poly contractible}
\po^{d,m}_n(\F)\cong \F^{md}
\quad
\mbox{ if }d<n.
\end{equation}
Because of this we only consider the case  
\begin{equation}\label{eq: d geq n}
d\geq n.
\end{equation}
}
\end{definition}
The following definitions will be used to formulate the notions of homotopy and homology stability.
\begin{definition}\label{def: basic definition}
{\rm
Let $f:X\to Y$ be a based map between based spaces $X$ and $Y$.
\par
(i)
Then the map $f$ 
is called {\it a homotopy equivalence}
(resp. {\it a homology equivalence})  
{\it through dimension} $N$
if the induced homomorphism
\begin{equation}
f_*:\pi_k(X)\to\pi_k(Y)
\qquad
(\mbox{resp. }f_*:H_k(X;\Z) \to H_k(Y;\Z))
\end{equation}
is an isomorphism for any integer $k\leq N$.\footnote{%
This definition was used in several papers concerning homotopy or homology stability
(e.g. \cite{BHMM}, \cite{BHMM2}).
}
\par
(ii)
Similarly,
the map $f$ 
is called {\it a homotopy equivalence}
(resp. {\it a homology equivalence})  
{\it up to dimension} $N$
if the induced homomorphism
\begin{equation}
f_*:\pi_k(X)\to\pi_k(Y)
\quad
(\mbox{resp. }f_*:H_k(X;\Z) \to H_k(Y;\Z))
\end{equation}
is an isomorphism for any integer $k< N$
and an epimorphism for $k=N$.\footnote{%
This definition first appeared and was used in \cite{Se}.
}
\par
(iii)
Let $G$ be a group and $f:X\to Y$ be a
$G$-equivariant based map between $G$-spaces $X$ and $Y$.
Then the map $f$ is called
a {\it $G$-equivariant homotopy equivalence through dimension} $N$
(resp. a {\it $G$-equivariant homotopy equivalence up to dimension} $N$)
if the restriction map 
\begin{equation}
f^H=f\vert X^H:X^H\to Y^H
\end{equation}
is a homotopy equivalence through dimension $N$ 
(resp. a homotopy equivalence up to dimension $N$)
for every subgroup $H\subset G$, where
$W^H$ denotes the $H$-fixed subspace of a $G$-space $W$ given by
\begin{equation}
W^H=\{x\in W: h\cdot x=x\mbox{ for every }h\in H\}.
\end{equation}
}
\end{definition}
\par
\begin{remark}
{\rm
(i)
When $X\subset Y$ and the map $f:X\stackrel{\subset}{\longrightarrow} Y$ is an inclusion map, the map $f$
is a homotopy equivalence up to dimension $N$ if and only if the pair $(Y,X)$ is $N$-connected.
Moreover, the map $f$ is a homotopy equivalence through dimension $N$ if and only if
it is a homotopy equivalence up to dimension $N$ and it induces a monomorphism
$f_*:\pi_N(X)\to \pi_N(Y).$

(ii)
If we set $\Z_2=\{\pm 1\}$,
by identifying $S^2$ with $\C\cup \infty$, the space $S^2$ acquires a $\Z_2$-action
induced by complex conjugation on $\C$.
The space $\CP^N$ also has a natural $\Z_2$-action induced by complex conjugation.
Thus, we can consider the space $(\Omega_d^2\CP^N)^{\Z_2}$ of based $\Z_2$-equivariant maps from $S^2$ to $\CP^N$.
It follows from \cite[Lemma 3.3]{KY13} that there is a homotopy equivalence\footnote{%
Although the case $N=1$ is not treated in (\ref{eq: 1.7}),
an analogous method shows that there is a homotopy equivalence
$(\Omega_d^2\CP^1)^{\Z_2}\simeq \Omega^2S^{3}\times \Z$
\ for $N=1$.}
\begin{equation}\label{eq: 1.7}
(\Omega_d^2\CP^N)^{\Z_2}\simeq \Omega^2S^{2N+1}\times \Omega S^N
\quad
\mbox{ for $N\geq 2.$}
\end{equation}
For completeness, we give its proof in
\S \ref{section: proofs}
(see Lemma \ref{lmm: Lemma 3.3 in KY13}).
\qed
}
\end{remark}
\par\vspace{1mm}\par
Now recall
the already obtained results for the space $\po^{d,m}_n(\R)$.
First, consider the cases
$mn=2\Leftrightarrow (m,n)=(2,1)$ or $(1,2)$.
The case $(m,n)=(2,1)$ is described in the following theorem. 
\begin{thm}[\cite{Br}, \cite{Se}; the case $(m,n)=(2,1)$]
\label{thm: Segal (m,n)=(2,1)}
\par
We make the identification $S^2=\C\cup \infty$ and let
$(\Omega^2_d\CP^1)^{\Z_2}_j$ denote the space of 
base-point preserving maps $S^2\to \CP^1$ of degree $d$ 
which commute with complex conjugation and have
degree $j$ when restricted to the real axis $S^1=\R^1\cup\infty$.
For each such $j$, let $\po^{d,2}_{1,j}(\R)$ denote the inverse image of
$(\Omega^2_d\CP^1)^{\Z_2}_j$ under the natural map $i^{d,2}_{1,\R}:\po^{d,2}_1(\R)\to (\Omega^2_d\CP^1)^{\Z_2}$.
\par
$\I$
The space $\po^{d,2}_1(\R)$
consists of $(d+1)$ connected components
$$
\{\po^{d,2}_{1,j}(\R):j=d-2k,
0\leq k\leq d\}.
$$ 
\par
$\II$
If $j=d-2k$ and $0\leq k\leq d$, the natural inclusion map
$$
i^{d,2}_{1,j}:
\po^{d,2}_{1,j}(\R) \stackrel{}{\longrightarrow}
(\Omega^2_d\CP^1)^{\Z_2}_j\simeq \Omega^2S^3
$$
is a homotopy equivalence up to dimension
$\frac{1}{2}(d-|j|)$.
\qed
\end{thm}
\begin{remark}
{\rm 
(i)
The assertion (i) of Theorem \ref{thm: Segal (m,n)=(2,1)} follows from \cite{Br}. The assertion (ii) follows from \cite[Proposition 7.2]{Se}. Moreover, by \cite[Proposition 7.1]{Se} there is a homotopy equivalence $(\Omega^2_d\CP^1)^{\Z_2}_j\simeq \Omega^2 S^3$.
\par
(ii)
If $j=d-2k$ and $0\leq k\leq d$, by \cite[Proposition 7.1 (cf. page 46)]{Se} we obtain
a homotopy equivalence
$\po^{d,2}_{1,j}(\R)
\simeq
\po^{\frac{1}{2}(d-|j|),2}_1(\C)$.
\qed
}
\end{remark}
Let $\lfloor x\rfloor$ denote the integer part of a real number $x$ and
recall that in the case $(m,n)=(1,2)$ we know the following (see \cite[Theorem 9.3]{KY13} for details).
\begin{thm}[\cite{KY13}; the case $(m,n)=(1,2)$]\label{thm: the case (m,n)=(1,2)}
Let $d\geq 2$ and let
$\po^{d,1}_{2,j}(\R)$ denote the subspace of $\po^{d,1}_2(\R)$ consisting of all monic
polynomials $f(z)\in \po^{d,1}_2(\R)$ of degree $d$ of the form
\begin{equation}
f(z)=\Big(\prod_{k=1}^{d-2j}(z-x_k)\Big)\Big(\prod_{k=1}^j(z-a_k)(z-\overline{a_k})\Big)
\end{equation}
such that
$(\{x_k\}_{k=1}^{d-2j},\{a_k\}_{k=1}^j)\in C_{d-2j}(\R) \times C_j({\rm H}_+)$. 
Here,
${\rm H}_+$  denotes the upper half plane in $\C$ given by
\begin{equation}\label{eq: H+}
{\rm H}_+=\{\alpha\in \C:\mbox{\rm Im }(\alpha)>0\}
\end{equation}
and we denote by $C_k(X)$
the unordered configuration space of $k$ distinct  points of  $X$ defined by
(\ref{eq: CX}).
\par\vspace{1mm}\par
$\I$
The space $\po^{d,1}_2(\R)$ consists of 
$(\lfloor d/2\rfloor +1)$ connected components
$
\{\po^{d,1}_{2,j}(\R):0\leq j\leq \lfloor d/2\rfloor\},
$
and
there is a homotopy equivalence
$$
\po^{d,1}_{2,j}(\R)\simeq K({\rm Br}(j),1)\quad
\mbox{ for each }0\leq j\leq \lfloor d/2\rfloor,
$$
where 
${\rm Br}(j)$ is the Artin braid group on $j$ strings.
\par
$\II$
The restriction of the natural map 
$$
i^{d,1}_{2,\R;j}=i^{d,1}_{2,\R}\vert \Po^{d,1}_{2,j}(\R):\Po^{d,1}_{2,j}(\R)\to \Omega^2_j\CP^1\simeq
\Omega^2_{j}S^2
\simeq \Omega^2S^3
$$
is a homology equivalence up to dimension $\lfloor j/2\rfloor$ if $j\geq 3$, and
it is a homotopy equivalence through dimension $1$ if $j=2$.
\qed
\end{thm}
Finally, consider the space $\po^{d,m}_n(\R)$ for $mn\geq 3$.
We recall the following results from \cite{KY13}.
\begin{thm}[\cite{KY13}; the case $mn\geq 3$]\label{thm: KY13}
Let $m,n,d\geq 1$ be positive integers satisfying
$mn\geq 3$ and $d\geq n$, and let $D(d;m,n)$ denote the positive integer
given by
\begin{equation}\label{eq: D(d;m,n)}
D(d;m,n)=(mn-2)(\lfloor d/n\rfloor+1)-1.
\end{equation}
\par
$\I$
The natural map
$($defined by  $($\ref{eq: mapjR-1}$))$
$$
i^{d,m}_{n,\R}:\Po^{d,m}_n(\R)\to (\Omega^2_d\CP^{mn-1})^{\Z_2}
\simeq
\Omega^2S^{2mn-1}\times \Omega S^{mn-1}
$$
is a homotopy equivalence up to dimension
$D(d;m,n)$
if $mn\geq 4$, and  a homology equivalence through dimension
$D(d;m,n)$ if $mn=3$.
\par
$\II$
Let
$$
s^{d,m}_{n,\R}:\Po^{d,m}_n(\R)\to \po^{d+1,m}_n(\R)
$$
denote
the stabilization map
$($defined by $($\ref{eq:  stabilization K}$)$ for $\K=\R)$.
\par
\begin{enumerate}
\item[{\rm (a)}]
Let $mn\geq 4$.
Then the map $s^{d,m}_{n,\R}$ is a homotopy equivalence
if $\lfloor \frac{d}{n}\rfloor=\lfloor \frac{d+1}{n}\rfloor$, and is a homotopy equivalence up to dimension
$D(d;m,n)$ otherwise.
\item[{\rm (b)}]
Let $mn=3$.
Then the map $s^{d,m}_{n,\R}$ is a homology equivalence
if $\lfloor \frac{d}{n}\rfloor=\lfloor \frac{d+1}{n}\rfloor$, and is a homology equivalence through dimension
$D(d;m,n)=\lfloor \frac{d}{n}\rfloor$ otherwise.

\end{enumerate}
\par\vspace{1mm}\par
$\III$
There is a stable homotopy equivalence
\begin{equation}\label{eq: stable homotopy equiv}
\po^{d,m}_n(\R)
\simeq_s
\Big(\bigvee_{i=1}^{\lfloor d/n\rfloor}S^{(mn-2)i}\Big)
\vee
\Big(
\bigvee_{i\geq 0,j\geq 1, i+2j\leq \lfloor d/n \rfloor}
\Sigma^{(mn-2)(i+2j)}D_j\Big),
\end{equation}
where 
 $D_j$ denotes the equivariant half smash product
defined in
(\ref{eq: D_k}).
\end{thm}
\begin{remark}
{\rm
The assertions (i) and (iii) of Theorem \ref{thm: KY13} follow from
\cite[Theorem 1.3]{KY13}. 
The assertion (ii) follows from \cite[Theorem 6.2, Corollary 6.4]{KY13}.
\qed
}
\end{remark}
\paragraph{1.2 The main results.}
It follows from the above theorems and \cite{KY8} that homology stability
holds for the space $\po^{d,m}_n(\F)$ when $\F=\R$ or $\C$.  We also know
that homotopy stability holds for $\po^{d,m}_n(\C)$ if and only if
$(m,n)\ne (1,2)$ (\cite{Se0}, \cite{Se}, \cite{KY8}).  For $\F=\R$ the
following results were known (see Theorems
\ref{thm: Segal (m,n)=(2,1)}, \ref{thm: the case (m,n)=(1,2)} and
\ref{thm: KY13}):
\begin{enumerate}
\item[(a)]
If $(m,n)=(1,2)$, homotopy stability does not hold for
$\po^{d,m}_n(\R)$.
\item[(b)]
If $(m,n)\ne (1,2)$ and either $mn\geq 4$ or $(m,n)=(2,1)$,
homotopy stability holds for $\po^{d,m}_n(\R)$.
\end{enumerate}
Thus the two borderline cases left by our earlier work are
\[
mn=3\quad\Longleftrightarrow\quad (m,n)=(3,1)\ \hbox{ or }\ (1,3).
\]
They belong to two rather different classical strands of the subject.  The
case $(m,n)=(3,1)$ is a space of polynomial tuples without common roots, of the
type considered by Segal \cite{Se}, whereas the case $(m,n)=(1,3)$ is a space of single monic polynomials 
without $3$-fold roots
closely related to the discriminant-complement
problems studied by Arnold and Vassiliev \cite{Ar,Va}.

In the present paper we settle the Segal case $(m,n)=(3,1)$.
Our proof uses an auxiliary model for $\po^{d,3}_1(\R)$ carrying a natural
$S^1$-action, obtained by rotating two of the three polynomial coordinates.
For the Arnold--Vassiliev case $(m,n)=(1,3)$ we have proved homotopy
stability by a different argument. More precisely, for $d\geq 3$ the natural
map $i^{d,1}_{3,\R}$ and the stabilization map $s^{d,1}_{3,\R}$
are homotopy equivalences up to dimension $\lfloor d/3\rfloor$ in the
sense of Definition~\ref{def: basic definition}. We shall give the proof
in a separate paper.

\begin{remark}
{\rm
(i)
For $(m,n)=(3,1)$ the stability dimension in
(\ref{eq: D(d;m,n)}) is
\begin{equation}
D(d;3,1)=d.
\end{equation}

(ii)
If $d\geq 1$, there is an isomorphism
\begin{equation}\label{eq: 1.13}
\pi_1(\po^{d,3}_1(\R))\cong\Z.
\end{equation}
This follows from \cite[Lemma 6.3]{KY13}; see also Lemma
\ref{lmm: KY13 (m,n)=(3,1)} below.  We shall also use the homotopy
equivalence
\begin{equation}\label{eq: 1.14}
\Omega S^2\simeq \Omega S^3\times S^1,
\end{equation}
which follows from the Hopf fibration $S^1\to S^3\to S^2$.
\qed
}
\end{remark}

To pass from homology stability to homotopy stability, we control the
deck action on homology of universal covers. This also determines the
fundamental-group action on higher homotopy groups in the required range.
\begin{definition}
{\rm
A path-connected space $X$ is said to be {\it simple up to dimension} $N$
if its fundamental group $\pi_1(X)$ acts trivially on $\pi_k(X)$ for all
$k<N$.\footnote{This definition is due to G. Segal \cite[page 61]{Se}.}
The space $X$ is {\it simple} if this action is trivial for every $k\geq 1$.
Since the action of $\pi_1(X)$ on itself is conjugation, a space which is
simple up to dimension $N>1$ has abelian fundamental group.
}
\end{definition}

We can now state the main results of this article.
\begin{thm}\label{thm: KY14: homotopy simple}
The space $\po^{d,3}_1(\R)$ is simple if $d\equiv 1\pmod 2$, and is
simple up to dimension $d$ if $d\equiv 0\pmod 2$.
\end{thm}

The covering-homology argument, together with the ordinary homology
stability in Theorem \ref{thm: KY13}, also gives the following result
directly.
\begin{thm}\label{thm: KY14: (m,n)=(3,1)}
$\I$
The natural map
\[
i^{d,3}_{1,\R}:\Po^{d,3}_1(\R)\longrightarrow
(\Omega^2_d\CP^{2})^{\Z_2}
\simeq \Omega^2S^{5}\times \Omega S^{2}
\simeq \Omega^2S^{5}\times \Omega S^{3}\times S^1
\]
is a homotopy equivalence up to dimension $d$.
\par
$\II$
The stabilization map
\[
s^{d,3}_{1,\R}:\po^{d,3}_1(\R)\longrightarrow
\po^{d+1,3}_1(\R)
\]
is a homotopy equivalence up to dimension $d$.
\end{thm}

Let $\Z_2=\{\pm1\}$ act on the spaces $\po^{d,3}_1(\C)$ and $\Omega^2_d\CP^2$
by complex conjugation.  The induced actions on
these spaces
satisfy
\begin{equation}\label{eq: fixed points}
\po^{d,3}_1(\C)^{\Z_2}=\po^{d,3}_1(\R),\qquad
(s^{d,3}_{1,\C})^{\Z_2}=s^{d,3}_{1,\R},\qquad
(i^{d,3}_{1,\C})^{\Z_2}=i^{d,3}_{1,\R}.
\end{equation}
Combining Theorem \ref{thm: KY14: (m,n)=(3,1)} with the corresponding
complex stability results and the fixed-point criterion of James and Segal
\cite{JS}, we obtain the following equivariant form.

\begin{crl}\label{crl: KY14: equiv}
The maps
\[
\begin{cases}
i^{d,3}_{1,\C}:\Po^{d,3}_1(\C)\longrightarrow
\Omega^2_d\CP^{2}\simeq \Omega^2S^5,\\
s^{d,3}_{1,\C}:\po^{d,3}_1(\C)\longrightarrow
\po^{d+1,3}_1(\C)
\end{cases}
\]
are $\Z_2$-equivariant homotopy equivalences up to dimension $d$.
\qed
\end{crl}

\paragraph{1.3 Organization.}
In \S \ref{section Definitions} we recall the basic definitions and the
natural and stabilization maps used below.  In
\S \ref{section: (m,n)=(3,1)} we study the basic properties of the space $\po^{d,3}_1(\R)$ and 
investigate the nilpotent result for the homology of its universal covering space 
(Theorem \ref{thm: simple up to dimension}).
By using this we
prove that it
is simple for odd $d$ and simple up to dimension $d$ for even $d$
(Corollaries \ref{crl: Polyd,3,1 simple} and
\ref{crl: simple even case}).  In \S \ref{section: proofs} we prove the
main theorems (Theorems \ref{thm: KY14: homotopy simple} and
\ref{thm: KY14: (m,n)=(3,1)})
and the homotopy equivalence (\ref{eq: 1.7}) (Lemma
\ref{lmm: Lemma 3.3 in KY13}).  Finally, in
\S \ref{section: optimality} we show that the stability range for $d=2$ and $d=3$
is optimal
(see
Lemmas \ref{lmm: optimality d=2}, \ref{lmm: the case d=3} and
 Remark \ref{remark: best possible}).


\section{Basic notations and definitions}\label{section Definitions}
\paragraph{2.1  Basic definitions and notations.}
We first recall some notation and basic definitions from \cite{KY13} needed to state and understand our results.
\begin{definition}
{\rm
From now on, let $X$ and $Y$ be based, path-connected spaces.
\par
(i)
Let
$\Map(X,Y)$ (resp. $\Map^*(X,Y)$) denote the space
consisting of all continuous maps
(resp. base-point preserving continuous maps) from $X$ to $Y$
with the compact-open topology.
\par
(ii)
For each element $D\in \pi_0(\Map^*(X,Y))$, let
$\Map^*_D(X,Y)$ denote the path component of
$\Map^*(X,Y)$ which corresponds to $D$.
For each integer $d\in \Z=\pi_0(\Map^*(S^2,\CP^N))$,
let $\Omega^2_d\CP^N=\Map^*_d(S^2,\CP^N)$ denote the path component of
$\Omega^2\CP^N$ of based maps
from $S^2$ to $\CP^N$ of degree $d$.
\par
(iii)
Let $F(X,k)$ denote 
{\it the ordered configuration space} of distinct $k$ points of the base space $X$
given by
\begin{equation}\label{eq: FX}
F(X,k)= \{(x_1,\cdots ,x_k)\in X^k:x_i\not= x_j\mbox{ if }i\not=j\}.
\end{equation}
\par
(iv)
The symmetric group $S_k$ on $k$ letters acts freely on this space by
permuting the coordinates. Let $C_k(X)$ be {\it the unordered configuration
space} of distinct points of $X$
given by the orbit space
\begin{equation}\label{eq: CX}
C_k(X)=F(X,k)/S_k.
\end{equation}
\par
(v)
The group $S_k$ also acts  on
the $k$-fold smash product 
\begin{equation}
X^{\wedge k}
=X\wedge  \cdots \wedge X
\quad
\mbox{($k$-times) }
\end{equation}
by permuting the coordinates.
Define the
equivariant half smash product $D_k(X)$ by
\begin{equation}\label{eq: DX}
D_k(X)=F(\C,k)\ltimes_{S_k} X^{\wedge k}
=(F(\C,k)\times_{S_k}X^{\wedge k})/(F(\C,k)\times_{S_k}*),
\end{equation}
where $*$ denotes the base point of $X^{\wedge k}$.\footnote{%
For a group $G$ and based $G$-spaces $(X,*)$ and $(Y,*)$, let
$X\wedge_GY$ denote the quotient space  
$X\wedge_GY=(X\times_GY)/((X\times_G *)\cup (*\times_G Y)).$
With this notation, we can write $D_k(X)=F(\C,k)_+\wedge_{S_k}X^{\wedge k}$,
where
$F(\C,k)_+=F(\C,k)\cup \{*\}$ (disjoint union).
}
In particular, for $X=S^1$ we write
\begin{equation}\label{eq: D_k}
D_k=D_k(S^1).
\end{equation}
}
\end{definition}
\begin{remark}
{\rm
It follows from \cite[page 462]{May2} 
(cf. \cite{CMM}, \cite{Milgram})
that there are homotopy equivalences
\begin{equation}\label{eq: D1, D2}
D_1\simeq S^1 \quad \mbox{ and }\quad D_2\simeq \Sigma \RP^2.
\end{equation}
 \par
(ii)
Let $\{X_d\}_{d=1}^{\infty}$ be  a sequence of connected spaces and
let
\begin{equation}\label{eq: seq maps}
X_1\stackrel{f_1}{\longrightarrow}
X_2\stackrel{f_2}{\longrightarrow}
X_3\stackrel{f_3}{\longrightarrow}
\cdots\cdots \stackrel{}{\longrightarrow}
X_d
\stackrel{f_d}{\longrightarrow}
X_{d+1}\stackrel{f_{d+1}}{\longrightarrow}
\cdots
\end{equation}
be a sequence of based continuous maps
such that each map $f_d$ is a homotopy equivalence
(resp. homology equivalence)  up to dimension $n(d).$
Let $X_{\infty}$ denote the colimit  (or homotopy colimit)
$\dis 
X_{\infty}=\lim_{d\to\infty}X_d$
taken over the continuous maps $\{f_d\}$.
\par
We say that
{\it homotopy stability} (resp. {\it homology stability})
holds for the space $X_d$ (or the map $f_d$)
if it is possible to choose the positive integer $n(d)$ such that
$\dis \lim_{d\to\infty}n(d)=\infty$. 
In this situation we also say that
{\it homotopy stability} (resp. {\it homology stability})
holds for the space $X_d$ (or
the natural map
$\iota_d:X_d\to X_{\infty}$).
\qed
}
\end{remark}

\paragraph{2.2  Spaces of non-resultant systems.}

Let $\N$ be the set of all positive integers.
From now on, let $d\in \N$, let
$(m,n)\in \N^2\setminus\{(1,1)\}$, and let
$\F$ be a field with algebraic closure
$\overline{\F}$.

\begin{definition}
{\rm
\par
Let $\P_d(\F)$ denote the space of all monic polynomials with coefficients
in $\F$
$f(z)=z^d+a_1z^{d-1}+\cdots +a_{d-1}z+a_d\in \F [z]$ of degree $d$.
Note that there is a natural homeomorphism
$\P_d(\F)\cong \F^d$ given by
\begin{equation}
f(z)=z^d+\sum_{k=1}^da_kz^{d-k}\mapsto (a_1,\cdots ,a_d).
\end{equation}
}
\end{definition}
\begin{remark}
{\rm
Recall that the classical resultant of a system of polynomials vanishes if and only if they have a common solution in an algebraically closed field containing the coefficients. Systems that have no common roots are called \lq\lq non-resultant\rq\rq. 
For this reason, 
we call the space $\po^{d,m}_n(\F)$
{\it the space of non-resultant systems of bounded multiplicity with coefficients in $\F$} (where $n$ is the multiplicity bound). 
\qed
}
\end{remark}
\paragraph{2.3 Natural maps and stabilization maps.}
In this section we always assume that $mn\geq 3$.
Here we briefly recall from \cite{KY13} the definitions of several maps
needed to state our results.  Some require no auxiliary choices; we call
these the natural maps.  The stabilization maps require choices, but their
homotopy classes are independent of those choices.
We use $\K$ to denote either $\R$ or $\C$.\footnote{%
The space $\po^{d,m}_n(\F)$ can be defined for any field $\F$.
This space was defined and used by Farb and Wolfson \cite{FW}
for a finite field $\F=\F_q$.
In this paper we only consider $\F=\R$ or $\C$, so we use the symbol $\K$
when a statement is valid in both cases.
}
\begin{definition}
{\rm
Let $\K=\R$ or $\C$, and
let $\Z_2=\{\pm 1\}$ denote the (multiplicative) cyclic group of order $2$.
\par
(i)
For a monic polynomial $f(z)\in \P_d(\K)$, let
$F_n(f)=F_n(f)(z) \in \P_d(\K)^n$ denote the $n$-tuple 
of monic polynomials of degree $d$ given by
\begin{equation}\label{eq: Fn}
F_n(f)(z)=
(f(z),f(z)+f^{\prime}(z),f(z)+f^{\prime\prime}(z),
\cdots ,f(z)+f^{(n-1)}(z)).
\end{equation}
For each $\alpha\in\C$, a polynomial $f(z)\in\P_d(\K)$ is not
divisible by $(z-\alpha)^n$ in $\C[z]$ if and only if
$F_n(f)(\alpha)\ne\mathbf{0}_n$, where $\mathbf{0}_n=(0,\ldots,0)\in\C^n$.
\par\vspace{1mm}\par
(ii)
When $\K=\C$, by identifying $S^2=\C \cup\infty$
we define {\it a  natural map}
\begin{equation}\label{eq: mapjC}
i^{d,m}_{n,\C}:\Po^{d,m}_n(\C)\to \Omega^2_d\CP^{mn-1}
\simeq \Omega^2S^{2mn-1}
\quad \mbox{by}
\end{equation}
\begin{equation*}\label{eq: def natural mapC}
i^{d,m}_{n,\C}(f)(\alpha)
=
\begin{cases}
[F_n(f_1)(\alpha):F_n(f_2)(\alpha):\cdots :F_n(f_m)(\alpha)]
&
\mbox{if }\alpha\in \C
\\
[1:1:\cdots :1] & \mbox{if }\alpha =\infty
\end{cases}
\end{equation*}
for $f=(f_1(z),\cdots ,f_m(z))\in \Po^{d,m}_n(\C)$
and $\alpha\in \C \cup \infty =S^2$,
where
we choose the points $\infty$ and $*=[1:1:\cdots :1]$ as 
the base-points of $S^2$ and $\CP^{mn-1}$, respectively.
\par\vspace{1mm}\par
(iii)
We  regard the spaces $S^2=\C\cup \infty$ and $\CP^{mn-1}$ as
$\Z_2$-spaces with actions  induced by complex conjugation.
Let  $(\Omega^2_d\CP^{mn-1})^{\Z_2}$ denote the space 
consisting of 
all $\Z_2$-equivariant based maps of degree $d$
$f:(S^2,\infty )\to( \CP^{mn-1},*)$.
\par\vspace{1mm}\par
(iv)
Since
\[
\Po^{d,m}_n(\R)\subset \Po^{d,m}_n(\C)
\quad\hbox{and}\quad
i^{d,m}_{n,\C}(\Po^{d,m}_n(\R))\subset
(\Omega^2_d\CP^{mn-1})^{\Z_2},
\]
we can define {\it a natural map}
\begin{align}\label{eq: mapjR-1}
&i^{d,m}_{n,\R}:\Po^{d,m}_n(\R) \to (\Omega_d^{2}\CP^{mn-1})^{\Z_2}
\qquad
\mbox{as the restriction}
\\
\nonumber
i^{d,m}_{n,\R}&=i^{d,m}_{n,\C}\vert
\Po^{d,m}_n(\R):\Po^{d,m}_n(\R)\to
(\Omega^2_d\CP^{mn-1})^{\Z_2}.
\end{align}
}
\end{definition}
Next, recall the definition of the stabilization maps.
\begin{definition}\label{def: stabd}
{\rm 
\par
Let $mn\geq 3$ and let $\K=\R$ or $\C$ as before.
For each integer $d\geq 1$,
let $\{x_{d,i}:1\leq i\leq m\} \subset (d,d+1)$ be any fixed
real numbers such that 
$x_{d,i}\not= x_{d,k}$ if $i\not= k$
and set
$\C_d=\{\alpha \in \C:\mbox{Re }(\alpha )<d\}$.
Choose a homeomorphism $\phi_d:\C\stackrel{\cong}{\longrightarrow}\C_d$
satisfying the following condition:
\begin{enumerate}
\item[($\dagger$)]
$\phi_d (\R)=(-\infty,d)\times \{0\}$, $\phi_d({\rm H}_+)=(-\infty, d)\times (0,\infty)$,
and $\phi_d (\overline{\alpha})=\overline{\phi_d (\alpha)}$
for any $\alpha\in \C$, 
\end{enumerate}
where ${\rm H}_+$ denotes the upper half plane in $\C$ as in (\ref{eq: H+})
and
we identify $\C=\R^2$ in the usual way.
We define {\it the stabilization map}
\begin{align}\label{eq:  stabilization K}
s^{d,m}_{n,\K}:\Po^{d,m}_n(\K)&\to \Po^{d+1,m}_n(\K)
\quad
\mbox{ by}
\\
\nonumber
s^{d,m}_{n,\K}(f_1(z),\cdots ,f_m(z))
&=
\big(
(z-x_{d,1})\widetilde{\phi}_d(f_1),\cdots ,
(z-x_{d,m})\widetilde{\phi}_d(f_m)\big)
\end{align}
for $(f_1(z),\cdots ,f_m(z))\in\Po^{d,m}_n(\K)$, where
we set 
\begin{equation}
\widetilde{\phi}_d(f)=\prod_{k=1}^d(z-\phi_d (x_k))
\quad
\mbox{if $f=f(z)=\prod_{k=1}^d(z-x_k)\in \P_d(\K).$}
\end{equation}
\par
}
\end{definition}
\begin{remark}
{\rm
The right-hand side of
\begin{equation}\label{eq: equivariant stab}
s^{d,m}_{n,\R}=(s^{d,m}_{n,\C})^{\Z_2}
\end{equation}
denotes the restriction of $s^{d,m}_{n,\C}$ to the fixed-point spaces.
For definiteness we use the horizontal compression
$\phi_d(x+iy)=d-e^{-x}+iy$. An equivariant isotopy of this compression
or a continuous variation of the added points through admissible choices
gives a homotopy of stabilization maps. See also
\cite[Definition 3.11]{KY9} and \cite[Lemma 7.1]{KY15}.
\qed
}
\end{remark}
\section{The case $(m,n)=(3,1)$}\label{section: (m,n)=(3,1)}

In this section we describe some basic properties of
the space $\po^{d,3}_1(\R)$.
\paragraph{3.1 Basic results about $\po^{d,3}_1(\R)$.}
We first recall the facts from \cite{KY13} that will be used below.
\begin{lemma}[\cite{KY13}]\label{lmm: KY13 (m,n)=(3,1)}
$\I$
There is an isomorphism
\[
\pi_1(\po^{d,3}_1(\R))\cong\Z.
\]

$\II$
The stabilization map
\[
s^{d,3}_{1,\R}:\po^{d,3}_1(\R)\longrightarrow
\po^{d+1,3}_1(\R)
\]
is a homology equivalence through dimension $d$.

$\III$
The maps $s^{d,3}_{1,\R}$ and $i^{d,3}_{1,\R}$ induce isomorphisms
\[
\begin{cases}
(s^{d,3}_{1,\R})_*:\pi_1(\po^{d,3}_1(\R))
\stackrel{\cong}{\longrightarrow}
\pi_1(\po^{d+1,3}_1(\R))\cong\Z,\\
(i^{d,3}_{1,\R})_*:\pi_1(\po^{d,3}_1(\R))
\stackrel{\cong}{\longrightarrow}
\pi_1((\Omega^2_d\CP^2)^{\Z_2})
\cong\pi_1(\Omega^2S^5\times\Omega S^3\times S^1)\cong\Z.
\end{cases}
\]
\end{lemma}
\begin{proof}
Assertion (i) is \cite[Lemma 6.3]{KY13}.  Assertion (ii) follows from
part (ii) of Theorem \ref{thm: KY13}, since $D(d;3,1)=d$.  Assertion
(iii) follows from \cite[Corollary 8.1]{KY13}, together with
(\ref{eq: 1.7}), (\ref{eq: 1.13}) and (\ref{eq: 1.14}).
\end{proof}

\paragraph{3.2 The spaces $\po^{d,3}_1(\R)$ and 
$\pomod{d}$.
}
Since there is a homotopy equivalence
\begin{equation}\label{eq: conditions 31}
\po^{1,3}_1(\R)\cong \R^3\setminus
\{(a,a,a):a\in \R\}\simeq S^1,
\end{equation}
it suffices to consider the space $\po^{d,3}_1(\R)$ only for $d\geq 2$.
\par\vspace{2mm}\par
We would like to define a natural $S^1$-action on $\po^{d,3}_1(\R)$, but there seems to be no way to do so directly. 
For $d\geq 1$, we introduce an auxiliary model 
$\pomod{d}$ 
of $\po^{d,3}_1(\R)$, which is homeomorphic to
$\po^{d,3}_1(\R)$ and carries a natural $S^1$-action
(see Lemma \ref{lmm: homeo} and Definition \ref{def: poly2}).
\par

\begin{definition}\label{def: poly1}
{\rm
(i)
For each $d\geq 1$,
let $\pomod{d}$ denote the space of $3$-tuples
$(f_1(z),f_2(z),f_3(z))\in \R[z]^3$
of polynomials with real coefficients
satisfying the following two conditions:
\begin{enumerate}
\item[(\ref{def: poly1}.1)]
$\max \{\deg (f_2(z)),\deg (f_3(z))\}<d$ and $f_1(z)$ is a monic polynomial of degree $d$, where $\deg (g(z))$ denotes the degree of  $g(z)\in \R [z]$.
\item[(\ref{def: poly1}.2)]
The polynomials $\{f_1(z),f_2(z),f_3(z)\}$ have no common root, that is:
\begin{equation*}
(f_1(\alpha),f_2(\alpha),f_3(\alpha))\not= (0,0,0)={\bf 0}_3
\qquad
\mbox{ for any }\alpha \in \C.
\end{equation*}
\end{enumerate}
\par
(ii) 
Define the map
$\varphi_{d}:
\po^{d,3}_1(\R) 
\stackrel{}{\longrightarrow}
\pomod{d}$ 
by
\begin{align}
\varphi_{d}({\rm f})
&=
(f_1(z),f_2(z)-f_1(z),f_3(z)-f_1(z))
\end{align}
for ${\rm f}=(f_1(z),f_2(z),f_3(z))\in\po^{d,3}_1(\R).$
%
}
\end{definition}
\begin{lemma}\label{lmm: homeo}
The map
$\varphi_{d}:
\po^{d,3}_1(\R) 
\stackrel{\cong}{\longrightarrow}
\pomod{d}$ 
is a homeomorphism.
\end{lemma}
\begin{proof}
Define the map $\psi :\pomod{d}\to \po^{d,3}_1(\R)$  by
$$
\psi (f_1(z),f_2(z),f_3(z))= (f_1(z),f_2(z)+f_1(z),f_3(z)+f_1(z)).
$$
Since $\varphi_d\circ \psi=\mbox{id}$ and $\psi \circ \varphi_d=\mbox{id}$,
the map $\varphi_d$ is a homeomorphism.
\end{proof}
\begin{definition}\label{def: poly2}
{\rm
Let $d\geq 1$.
\par
(i)
 Define an
$S^1$-action on the space $\pomod{d}$ by
\begin{equation}\label{eq: action def}
e^{\sqrt{-1}\ \theta}\cdot
{\rm f}
=
(f_1(z),g(z),h(z))
\end{equation}
for $\theta\in \R$ and
${\rm f}=(f_1(z),f_2(z),f_3(z))\in \pomod{d}$, where the real coefficients 
polynomials $g(z)$ and $h(z)\in \R[z]$ are given by
\begin{align}
\begin{pmatrix}
g(z)
\\
h(z)
\end{pmatrix}
&=
\begin{pmatrix}
\cos \theta & -\sin \theta
\\
\sin \theta & \cos \theta
\end{pmatrix}
\begin{pmatrix}
f_2(z)
\\
f_3(z)
\end{pmatrix}
=
\begin{pmatrix}
f_2(z)\cos \theta-f_3(z)\sin \theta
\\
f_2(z)\sin \theta+f_3(z)\cos \theta
\end{pmatrix}
\\
\nonumber
&\Leftrightarrow
\ \ 
g(z)+\sqrt{-1}\ h(z)=e^{\sqrt{-1}\theta}\big(f_2(z)+\sqrt{-1}\ f_3(z)\big).
\end{align}
\par
(ii)
The $S^1$-action on $\pomod{d}$ is free: a nonidentity plane rotation can fix the pair $(f_2,f_3)$ only when both polynomials are identically zero, which is excluded by the common-root condition. We use its homotopy orbit space.
Let 
$(\pomod{d})_{S^1}$ be the space defined by the Borel construction 
\begin{equation}
(\pomod{d})_{S^1}=ES^1\times_{S^1}\pomod{d}.
\end{equation}
}
\end{definition}
\begin{example}
{\rm 
Let $d\geq 2$, and set ${\rm f}_0=(z^d,1,z)\in \pomod{d}$. 
Then we have:
\begin{align}\label{eq: example: action of S1}
e^{\sqrt{-1}\theta}\cdot {\rm f}_0&=(z^d,\cos \theta -z\sin \theta,\sin \theta +z\cos \theta)
=(z^d,g_{\theta}(z),h_{\theta}(z)),
\end{align}
where we set 
$(g_{\theta}(z),h_{\theta}(z))=(\cos \theta -z\sin \theta ,\sin \theta +z\cos \theta).$
We then have
\begin{equation}\label{eq: g+h}
e^{\sqrt{-1}\theta}=
g_{\theta}(0)+\sqrt{-1}h_{\theta}(0).
\qquad
\quad \qed
\end{equation}
}
\end{example}
Since $S^1$ acts on the space $\pomod{d}$, 
we obtain a fibration sequence
\begin{equation}\label{eq: fibration qd}
S^1 \stackrel{\hat{i}_d}{\longrightarrow}
\pomod{d}
\stackrel{\hat{q}_d}{\longrightarrow}
(\pomod{d})_{S^1},
\end{equation}
where 
$\hat{q}_d$ denotes the natural projection and
$\hat{i}_d$ is the fiber inclusion represented by the orbit of ${\rm f}_0$
as in (\ref{eq: example: action of S1}), i.e.
\begin{equation}
\hat{i}_d(e^{\sqrt{-1}\theta})=(z^d,g_{\theta}(z),h_{\theta}(z))\qquad(d\geq2)
\quad\quad
\mbox{for }\theta \in \R.
\end{equation}
For $d=1$ use the orbit of $(z,1,0)$, given by $(z,\cos\theta,\sin\theta)$.
\begin{definition}
{\rm
Let $d\geq 1$ be a positive integer.
\par 
(i)
First, when $d\equiv 1$ $(\mbox{mod }2)$,
define a map $\rho_d:\pomod{d}\to S^1$ by
\begin{equation}\label{def: rho d odd}
\rho_d({\rm f})=
\prod_{j=1}^l
\Bigg(
\frac{
f_2(x_j)+\sqrt{-1}f_3(x_j)}
{\vert f_2(x_j)+\sqrt{-1}f_3(x_j)\vert}
\Bigg)^{\epsilon (j)}
\end{equation}
for  ${\rm f}=(f_1(z),f_2(z),f_3(z))\in \pomod{d}$, where $\epsilon (j)=(-1)^{j-1}$ and
the monic polynomial $f_1(z) \in \R[z]$ is represented in the following form:
\begin{enumerate}
\item[$(\dagger\dagger)$]
There exists some real numbers
$x_1\leq x_2\leq \cdots \leq x_l$ such that
\begin{equation}
f_1(z)=(z-x_1)(z-x_2)\cdots (z-x_l)g(z)
\end{equation}
and
$g(z)\in \R [z]$ is a monic polynomial of degree $d-l$ without a real root.
\end{enumerate}
\par
If $x_j=x_{j+1}$, then
$$
\big(f_2(x_j)+\sqrt{-1}f_3(x_j)\big)^{\epsilon (j)}
\big(f_2(x_{j+1})+\sqrt{-1}f_3(x_{j+1})\big)^{\epsilon (j+1)}=1.
$$
Moreover, since $d\equiv 1$ $(\mbox{mod }2)$,
the polynomial $f_1(z)$ always has a real root, and
the polynomial $f_1(z)$ can be represented  in the form $(\dagger\dagger)$.
Thus,
the map $\rho_d$ is well defined. It is also continuous:
near a real root of multiplicity $m$, the number of nearby real roots
counted with multiplicity has parity $m$, since nonreal roots occur
in conjugate pairs. The evaluation phases approach the same nonzero
phase. Adjacent equal limiting phases cancel in the alternating
product, and insertion or removal of two real roots does not change
the signs of subsequent factors. The same argument applies in even degree.
\par
(ii)
When $d\equiv 0$ $(\mbox{mod }2)$,
define the map  $\rho_d:\pomod{d}\to S^1$ by
\begin{equation}\label{def: rho d even}
\rho_d({\rm f})=
\begin{cases}
\dis
\prod_{j=1}^l
\Bigg(
\frac{
f_2(x_j)+\sqrt{-1}f_3(x_j)}
{\vert f_2(x_j)+\sqrt{-1}f_3(x_j)\vert}
\Bigg)^{\epsilon (j)}
& \mbox{if  $(\dagger\dagger)$ is satisfied,}
\\
1 & \mbox{otherwise}
\end{cases}
\end{equation}
for  ${\rm f}=(f_1(z),f_2(z),f_3(z))\in \pomod{d}$.
\par
The same continuity argument applies when $d$ is even.

\par
(iii)
We also define two maps 
\begin{equation}
\begin{cases}
q_d:\po^{d,3}_1(\R)\to (\pomod{d})_{S^1}
\\
r_d:\po^{d,3}_1(\R)\to S^1
\end{cases}
\qquad
\mbox{by}
\end{equation}
\begin{equation}\label{eq: varphi}
q_d=\hat{q}_d\circ \varphi_d
\  \ \mbox{ and }\ 
\quad
r_d=\rho_d\circ \varphi_d.
\end{equation}
\par
(iv)
Given two maps
$
Y \stackrel{f}{\longleftarrow} X\stackrel{g}{\longrightarrow} Z,
$
let 
$(f,g):X\to Y\times Z$ denote the map defined by
\begin{equation}
(f,g)(x)=(f(x),g(x))
\qquad
\mbox{ for }x\in X.
\end{equation}
\par
(v) Let 
\begin{equation}
u_d:\widetilde{\po}^{d,3}_1\to \po^{d,3}_1(\R)
\end{equation}
denote the universal covering of the space $\po^{d,3}_1(\R)$.
}
\end{definition}
\begin{lemma}\label{lmm: polyd: odd}
Let $d\geq 3$ such that
$d\equiv 1\ \mbox{\rm (mod }2)$.
\par\vspace{1mm}\par
$\I$ 
The space $(\pomod{d})_{S^1}$ is simply connected, and
the map
\begin{equation}\label{eq: homotopy decomp.}
(q_d,r_d):
\po^{d,3}_1(\R)\stackrel{\simeq}{\longrightarrow} (\pomod{d})_{S^1}\times S^1
\end{equation}
is a homotopy equivalence.
\par
$\II$
The homomorphism
$$
(q_d)_*:\pi_k(\po^{d,3}_1(\R))
\stackrel{\cong}{\longrightarrow}
\pi_k((\pomod{d})_{S^1})
$$
is an isomorphism for any $k\geq 2$.
\par
$\III$
The map
$$
q_{d}\circ u_d:
\widetilde{\po}^{d,3}_1
\stackrel{\simeq}{\longrightarrow}
(\pomod{d})_{S^1}
$$
is a homotopy equivalence.
\par
$\IV$
There is a homotopy equivalence
\begin{equation}\label{eq: univ. cover for d:odd}
\po^{d,3}_1(\R)\simeq \widetilde{\po}^{d,3}_1\times S^1.
\end{equation}
\end{lemma}
\begin{proof}
(i)
Using (\ref{eq: g+h}), we obtain
\begin{equation}\label{eq: Rdjd}
\rho_d\circ \hat{i}_d=\mbox{id}.
\end{equation}
Since $\pi_1(\pomod{d})\cong \pi_1(\po^{d,3}_1(\R))\cong \Z$,
from the long exact homotopy sequence of the fibration
(\ref{eq: fibration qd}),
we obtain the following equality:
\[
(\rho_d)_*\circ (\hat{i}_d)_*=\mathrm{id}:
\Z=\pi_1(S^1) \xrightarrow{(\hat{i}_d)_*} \pi_1(\pomod{d})=\Z \xrightarrow{(\rho_d)_*} \pi_1(S^1)=\Z.
\]
Thus, it follows from (\ref{eq: fibration qd}) that we have
the following three assertions: 
\begin{enumerate}
\item[(a)]
The  two homomorphisms
\begin{equation}\label{eq: rho pi1}
\begin{cases}
(\rho_d)_*:\pi_1(\pomod{d})\stackrel{\cong}{\longrightarrow}\pi_1(S^1)
\\
(\hat{i}_d)_*:\pi_1(S^1)\stackrel{\cong}{\longrightarrow}\pi_1(\pomod{d})
\end{cases}
\end{equation}
are isomorphisms.
\item[(b)]
The space $(\pomod{d})_{S^1}$ is simply connected.
\item[(c)]
The homomorphism
$(\hat{q}_d)_*:\pi_k(\pomod{d})\stackrel{\cong}{\longrightarrow} \pi_k((\pomod{d})_{S^1})$
is an isomorphism for any $k\geq 2$.
\end{enumerate}
It follows that the map $(\hat{q}_d,\rho_d)$ induces an isomorphism
\begin{equation}\label{eq; homotopy equiv. univ}
(\hat{q}_d,\rho_d)_*:
\pi_k(\pomod{d})
\stackrel{\cong}{\longrightarrow}
\pi_k((\pomod{d})_{S^1}\times S^1)
\end{equation}
for every $k$. Thus the map
$(\hat{q}_d,\rho_d)$
is a homotopy equivalence.
Therefore, by using the homeomorphism $\varphi_d$
(by Lemma \ref{lmm: homeo}), we also obtain a homotopy equivalence
\begin{equation*}\label{eq: poly homo equiv}
\begin{CD}
(\hat{q}_d,\rho_d)\circ \varphi_d:
\po^{d,3}_1(\R) @>\varphi_d>\cong> \pomod{d}
@>(\hat{q}_d,\rho_d)>\simeq>
(\pomod{d})_{S^1}\times S^1.
\end{CD}
\end{equation*}
Since $(\hat{q}_d,\rho_d)\circ \varphi_d=
(\hat{q}_d\circ \varphi_d,\rho_d\circ \varphi_d)=(q_d,r_d)$,
the assertion (i) is proved.
\par
(ii)
The assertion (ii) follows from (c) and (\ref{eq: varphi}).
\par
(iii)
Consider the composite of maps
$$
q_{d}\circ u_d:
\widetilde{\po}^{d,3}_1
\stackrel{u_d}{\longrightarrow}
\po^{d,3}_1(\R)
\stackrel{q_d}{\longrightarrow}
(\pomod{d})_{S^1}.
$$
Since $u_d$ is a covering projection,
it induces an isomorphism on the homotopy group
$\pi_k(\ )$ for any $k\geq 2$.
Thus, by (ii),
the map $q_d\circ u_d$ induces an isomorphism on $\pi_k(\ )$ for any $k\geq 2$.
Since both spaces $\widetilde{\po}^{d,3}_1$ and $(\pomod{d})_{S^1}$ are simply connected, 
the map $q_d\circ u_d$ is indeed a homotopy equivalence.
\par
(iv) The assertion easily follows from the assertions (i) and (iii).
\end{proof}
\begin{crl}\label{crl: Polyd,3,1 simple}
If $d\equiv 1\ \mbox{\rm (mod }2)$,
the space $\po^{d,3}_1(\R)$ is simple.
\end{crl}
\begin{proof}
Since $S^1$ is a Lie group, it is simple.
Then the assertion follows easily from the homotopy equivalence (\ref{eq: univ. cover for d:odd}).
\end{proof}
\begin{lemma}\label{lmm: 3.9}
$\I$ For each $\theta \in \R$ and ${\rm f}=(f_1(z),f_2(z),f_3(z))\in \pomod{d}$,
$$
\rho_d(e^{\sqrt{-1}\theta}\cdot {\rm f})=
\begin{cases}
e^{\sqrt{-1}\theta}\rho_d({\rm f}) & \mbox{ if }d\equiv 1 \ (\mbox{{\rm mod }$2)$}
\\
\rho_d({\rm f}) & \mbox{ if }d\equiv 0 \ (\mbox{{\rm mod }$2)$}
\end{cases}
$$

\par
$\II$
The induced homorphism
$$
(\rho_d)_*:\pi_1(\pomod{d})
\stackrel{\cong}{\longrightarrow}
\pi_1(S^1)
$$
is an isomorphism.
\end{lemma}
\begin{proof}
(i) If $f_1(z)$ has a complex root $\alpha\in \C\setminus \R$,
its complex conjugate $\overline{\alpha}$ is also a root of $f_1(z)$.
Thus, the number of real roots of $f_1(z)$, counted with multiplicity, is odd if $d\equiv 1$ (mod $2$) and it is even if $d\equiv 0$ (mod $2$).
Thus, the assertion (i) easily follows from the definition of $\rho_d$.
\par
(ii)
Note that  the assertion (ii) was already proved if $d\equiv 1$ (mod $2$) in (\ref{eq: rho pi1}).
Now we give another easier proof, which works for the case $d\equiv 0$ (mod $2$), too.
Let $x_1<x_2<\cdots < x_d$ be any fixed $d$ distinct real numbers and let
$f_1(z)=\prod_{k=1}^d (z-x_k)\in \R[z]$.
For each $\theta \in \R$,
by using real Lagrange interpolation we choose the real coefficients polynomials $f_{2,\theta}(z)$ and $f_{3,\theta}(z)$ 
of degree $<d$ such that
\begin{equation}\label{eq: condtion gmh}
(f_{2,\theta}(x_k),f_{3,\theta}(x_k))=
\begin{cases}
(\cos \theta ,\sin \theta) & \mbox{ if }k=1
\\
(1,0) &\mbox{ if }2\leq k\leq d
\end{cases}
\end{equation}
More explicitly, put
\[
L_1(z)=\prod_{k=2}^d\frac{z-x_k}{x_1-x_k}.
\]
The polynomials
\[
f_{2,\theta}(z)=1+(\cos\theta-1)L_1(z),\qquad
f_{3,\theta}(z)=\sin\theta\,L_1(z)
\]
satisfy (\ref{eq: condtion gmh}) and depend continuously on $\theta$.
\par
Let ${\rm f}_{\theta}=(f_1(z),f_{2,\theta}(z),f_{3,\theta}(z))\in \pomod{d}$.
Then by (\ref{eq: condtion gmh}),
$$
\rho_d ({\rm f}_{\theta})
=\prod_{k=1}^d
\Bigg(
\frac{
f_{2,\theta}(x_k)+\sqrt{-1}\ f_{3,\theta}(x_k)
}
{
\vert f_{2,\theta}(x_k)+\sqrt{-1}\ f_{3,\theta}(x_k)\vert
}\Bigg)^{\epsilon (k)}=e^{\sqrt{-1}\theta}.
$$
Thus, $(\rho_d)_*:\pi_1(\pomod{d})\to \pi_1(S^1)=\Z$
is an epimorphism.
However,
since $\pi_1(\pomod{d})=\Z$, $(\rho_d)_*$ is an isomorphism.
\end{proof}
\begin{remark}
$\I$
If $d\equiv 1$ $(${\rm mod }$2)$, there is a homeomorphism
\begin{equation}\label{eq: 3.23}
\pomod{d} \cong \rho_d^{-1}(1)\times S^1
\quad
\mbox{given by }\ \ 
{\rm f}\mapsto (\rho_d({\rm f})^{-1}\cdot {\rm f},\rho_d({\rm f})),
\end{equation}
 and there is a homotopy equivalence
\begin{equation}\label{eq: 3.24}
\widetilde{\po}^{d,3}_1\simeq \rho_d^{-1}(1).
\end{equation} 
\par
$\II$
Let $d\equiv 0$ $(${\rm mod }$2)$, and let
\begin{equation}\label{eq: seq fib}
S^1 \stackrel{i}{\longrightarrow} 
\pomod{d}
\stackrel{q_d}{\longrightarrow}
(\pomod{d})_{S^1}
\end{equation}
denote the fibration sequence.
\par
Then the induced homorphism
$
i_*:\pi_1(S^1)\to \pi_1(\pomod{d})
$
is trivial, and
$\pi_1((\pomod{d})_{S^1})=\Z$.
Hence, the space $(\pomod{d})_{S^1}$  is not the universal covering of $\po^{d,3}_1(\R)$.
\end{remark}
\begin{proof}
(i)
The first assertion can be easily seen by using (i) of Lemma \ref{lmm: 3.9}.
Since $\pi_1(\pomod{d})=\Z$, we also obtain the second assertion.
\par
(ii)
Next suppose that $d=2d_1\equiv 0$ (mod $2$), and consider the fibration sequence
(\ref{eq: seq fib}).
Let  $f_0(z)=(z^2+1)^{d/2}=(z^2+1)^{d_1}$, and let
 $(g(z),h(z))\in \R[z]^2$ be the pair of polynomials such that
 $\max\{\deg g, \deg h\}<d$.
Since
$$
(f_0,g,h)\in \pomod{d}
\Leftrightarrow (g(\sqrt{-1}),h(\sqrt{-1}))\in \C^2\setminus \{0\},
$$
the allowed fixed $f_0$-fiber
$F(f_0)=
\{(g,h)\in \R [z]^2:(f_0,g,h)\in \pomod{d}\}
$
is homeomorphic to the space
$(\C^2\setminus \{0\})\times \R^{2d-4}$.
Thus, it is simply connected and contains the orbit of $(g,h)=(1,0)$.
Hence, this orbit contracts.
Since any orbit classes at other base-points are related by the homotopy obtained from a path and the $S^1$-action,
all are null-homotopic.
Hence, by using the homotopy exact sequence induced from (\ref{eq: seq fib}), we see that
$i_*$ is trivial, and that
$(q_d)_*:\Z=\pi_1(\pomod{d})\stackrel{\cong}{\longrightarrow}
\pi_1((\pomod{d})_{S^1})$
is an isomorphism.
Thus, (ii) was proved.
\end{proof}

\paragraph{3.3 \ Nilpotent result.}
We shall use a finite degree filtration to control the deck action on
the homology of the universal covering. 
\par
\begin{example}
{\rm
Let $X$ be a connected CW complex with $\pi_1(X)=\Z$ and let $\widetilde{X}\to X$
be a universal covering of the space $X$.
If we identify $\pi_1=\pi_1(X)=\Z$ with the infinite cyclic group $\langle t\rangle$ generated by the element $t$,
one can also identify its group ring $\Z [\pi_1]$ with the Laurent polynomial ring
\begin{equation}
\Lambda=
\Z[\pi_1]=\Z[t,t^{-1}].
\end{equation}
Note that the homology group $H_*(\widetilde{X};\Z)$ is a $\Z[t,t^{-1}]$-module.
\qed
}
\end{example}
Recall the following definition.

\begin{definition}
{\rm 
A module $M$ over
$\Lambda=\Z[t,t^{-1}]$ will be called {\it nilpotent} if 
$(t-1)^aM=0$ for some integer $a\geq 1$. Submodules, quotients and extensions of nilpotent modules are nilpotent.
%
}
\end{definition}
\begin{thm}\label{thm: simple up to dimension}
$\I$
If $d\equiv 1$ $(${\rm mod }$2)$, 
the deck action of
$\pi_1(\po^{d,3}_1(\R))\cong\Z$ on the homology 
$H_k(\widetilde{\po}^{d,3}_1;\Z)$ is trivial for any $k\geq 0$.
\par
$\II$ 
If $d\equiv 0$ $(${\rm mod }$2)$, 
the deck action of
$\pi_1(\po^{d,3}_1(\R))\cong\Z$ on the homology
$H_k(\widetilde{\po}^{d,3}_1;\Z)$ is nilpotent for any $0\leq k<d$.
%
%
\end{thm}

\paragraph{3.4 Preliminary results for the proof of Theorem \ref{thm: simple up to dimension}.}
\begin{remark}\label{remark :3.14}
{\rm 
(i)
We can easily prove (i) of Theorem \ref {thm: simple up to dimension}  
by using
Corollary \ref{crl: Polyd,3,1 simple}.
So it remains to prove the case $d\equiv 0$ (mod $2$) with $d\geq 2$.
\par
(ii)
For $(f,g,h)\in \pomod{d}$, we frequently identify the real polynomial pair $(g,h)\in \R[z]^2$ with the
complex polynomial 
\begin{equation}\label{eq: q-h+ih}
q=q(z)=g+\sqrt{-1}h\in \C[z].
\end{equation}
\par
(iii)
Remark that $(g,h)\not= (0,0)$ if
$(f,g,h)\in \pomod{d}.$
Thus, if $(f,g,h)\in \pomod{d}$ and $q=q(z)=g+\sqrt{-1}h$, 
then $0\leq \deg (q)\leq d-1$.
\qed
}
\end{remark}
\begin{definition}
{\rm 
Put $F_{-1}=\emptyset$. For $0\leq r\leq d-1$, define $F_r$ and $S_r$ by
\begin{equation}
\begin{cases}
F_r&=\{(f,g,h)\in \pomod{d}:q=g+\sqrt{-1}h,\deg (q)\leq r\},
\\
S_r&= F_r\setminus F_{r-1}.
\end{cases}
\end{equation}
Note that there is the following filtration
\begin{equation}
\emptyset=F_{-1}\subset F_0\subset\cdots\subset F_{d-1}=\pomod{d}.
\end{equation}
}
\end{definition}
\begin{remark}
{\rm 
(i)
Each $F_r$ is open in the linear coefficient space 
\begin{equation}
L_r=\R^d\times \C^{r+1}
\end{equation}
 and closed in $F_{r+1}$.
In particular, $F_r$ is an open smooth manifold  of dimension $d+2(r+1)$.
\par
(ii)
For each $1\leq r\leq d-1$, the submanifold $F_{r-1}\subset F_r$ has real codimension $2$ 
with normal framing  given by the two real components of coefficient $z^r$ in $q=g+\sqrt{-1}h$.
\par
(iii)
Let $1\leq r\leq d-1$ and $(f,g,h)\in S_r$ and we write
$q=q(z)=g+\sqrt{-1}h\in \C[z]$ as in (\ref{eq: q-h+ih}).
Then we can uniquely write
\begin{equation}\label{eq: G H}
q=g+\sqrt{-1}h=u(G+\sqrt{-1}H)
\qquad
(u\in \C^*),
\end{equation}
where
$G=G(z)\in \R[z]$ is a monic polynomial of $\deg (G)=r$ and $H=H(z)\in \R[z]$
is a polynomial of $\deg (H)<r$.
\par
By using polynomial division, we can also write
\begin{equation}\label{eq: f=AG+R}
f=AG+R,
\end{equation}
where $A=A(z)\in \R[z]$ is monic polynomial  of $\deg (A)=d-r$ and 
$R=R(z)\in \R[z]$ of $\deg (R)<r$.
\par
(iv)
Multiplication by $u\in\C^*$ acts on the pair $(G,H)$ by an invertible
real $2\times2$ matrix. This matrix remains invertible on evaluating
the polynomials at any complex number. Hence
\begin{equation}\label{eq: equiv}
g(\alpha)=h(\alpha)=0
\ \ 
\Leftrightarrow
\ \ 
G(\alpha)=H(\alpha)=0 \quad
\end{equation}
for $\alpha\in \C$.
At a common zero of $G$ and $H$ we have $f(\alpha)=R(\alpha)$, so
\begin{equation}\label{eq: (G,H,R)}
(G,H,R)\in \pomod{r}.
\end{equation}
By using this correspondence $(f,g,h)\mapsto (u,A,(G,H,R))$,
we finally obtain the homeomorphism
\begin{equation}\label{eq: homeo-Sr}
S_r\cong \C^*\times \R^{d-r}\times \pomod{r}
\quad
\mbox{ for $1\leq r\leq d-1$.}
\end{equation}
In particular, $S_0\cong \C^*\times \R^d$
and
the real codimension of $S_r$ in $\pomod{d}$ is
\begin{equation}
c_r=2(d-1-r).
\end{equation}
\par
(v)
By (\ref{eq: G H}), (\ref{eq: f=AG+R}) and (\ref{eq: (G,H,R)}), one can define the map
\begin{equation}\label{PhiA}
\Phi_A:\pomod{r} \to \pomod{d}
\quad \mbox{ by }\quad
\Phi_A
(G,H,R)= (AG+R,G,H).
\qed
\end{equation}
}
\end{remark}
\begin{lemma}\label{lmm: pi1 induced}
If $1\leq r\leq d-1$, the map $\Phi_A$ of (\ref{PhiA}) induces an isomorphism
\[
(\Phi_A)_*:\pi_1(\pomod{r})\stackrel{\cong}{\longrightarrow}
\pi_1(\pomod{d})\cong\Z.
\]
\end{lemma}
The map $\Phi_A$ is the slice $u=1$ of the parametrization
\eqref{eq: homeo-Sr}, with $A$ fixed.

We first record the elementary root-perturbation calculation needed
in the proof. Choose $d$ distinct real numbers $0=a_1<a_2<\cdots<a_d$
and put
\[
G(z)=\prod_{k=1}^r(z-a_k),\qquad
A(z)=\prod_{k=r+1}^d(z-a_k).
\]
Let $H_0$ be the polynomial of degree less than $r$ satisfying
$H_0(a_1)=0$ and $H_0(a_k)=1$ for $2\leq k\leq r$, and put
$L=1-H_0$. In particular, $H_0=0$ and $L=1$ if $r=1$.
For $\theta\in\R$ and $\epsilon>0$, set
\[
\lambda=\epsilon\cos\theta,\qquad \mu=\epsilon\sin\theta,\qquad
H_\theta=H_0+\lambda L.
\]
For all sufficiently small $\epsilon\geq0$, uniformly in $\theta$,
$AG+\mu L$ has $d$ distinct real roots
$\alpha_1(\epsilon,\theta)<\cdots<\alpha_d(\epsilon,\theta)$
with $\alpha_k(0,\theta)=a_k$. We suppress $\theta$ in the notation
for these roots.
\begin{lemma}\label{lmm: estimation of G, H}
For $1\leq k\leq d$,
\[
\alpha_k(\epsilon)=a_k-\frac{\epsilon\sin\theta\,L(a_k)}
{(AG)'(a_k)}+O(\epsilon^2).
\]
In particular, if $k=1$,
\[
\begin{pmatrix}G(\alpha_1(\epsilon))\\ H_\theta(\alpha_1(\epsilon))\end{pmatrix}
=
\begin{pmatrix}0&-1/A(0)\\1&-H_0'(0)/(A(0)G'(0))\end{pmatrix}
\begin{pmatrix}\epsilon\cos\theta\\\epsilon\sin\theta\end{pmatrix}
+O(\epsilon^2).
\]
For $k>1$, the pair $(G(\alpha_k(\epsilon)),H_\theta(\alpha_k(\epsilon)))$
converges uniformly in $\theta$ to the nonzero pair $(G(a_k),H_0(a_k))$.
\end{lemma}
\begin{proof}
Apply the implicit function theorem at $(a_k,0)$ to
\[
F(z,\epsilon)=A(z)G(z)+\epsilon\sin\theta\,L(z).
\]
This gives
\[
\alpha_k'(0)=-\frac{L(a_k)\sin\theta}{(AG)'(a_k)}.
\]
The denominator is nonzero at each $a_k$, since all roots of $AG$
are simple. Taylor expansion gives the stated formula, uniformly
in $\theta$ by compactness of the parameter circle. For $k=1$ we have
$L(0)=1$, $G(0)=H_0(0)=0$ and $(AG)'(0)=A(0)G'(0)$, so substitution
gives the displayed matrix. For $2\leq k\leq r$ the limiting pair
is $(0,1)$, whereas for $k>r$ its first coordinate $G(a_k)$ is nonzero.
\end{proof}
\begin{proof}[Proof of Lemma \ref{lmm: pi1 induced}]
The space of monic polynomials $A$ of degree $d-r$ is contractible,
so we may use the choice above. The loop
\[
\omega(\theta)=(G,H_\theta,\mu L)\in\pomod{r}
\]
has evaluation $H_\theta(a_k)+\sqrt{-1}\mu L(a_k)$ equal to
$\epsilon e^{\sqrt{-1}\theta}$ for $k=1$ and to $1$ for
$2\leq k\leq r$. Thus $\rho_r(\omega(\theta))=e^{\sqrt{-1}\theta}$,
and Lemma \ref{lmm: 3.9} shows that $\omega$ represents a generator
of the fundamental group $\pi_1(\pomod{r})\cong \Z$.
Its image is
\[
\Phi_A(\omega(\theta))=(AG+\mu L,G,H_\theta).
\]
The matrix in Lemma \ref{lmm: estimation of G, H} has determinant
$1/A(0)\ne0$. For sufficiently small $\epsilon$ the evaluation loop
at the first continued root therefore has winding number $\pm1$,
by homotopy to the corresponding ellipse. Every other evaluation
loop lies near a fixed nonzero pair and has winding number zero.
By the alternating-product definition of $\rho_d$, the loop
$\rho_d\circ\Phi_A\circ\omega$ has degree $\pm1$. Since
$(\rho_d)_*$ is an isomorphism by Lemma \ref{lmm: 3.9}, this proves
the assertion.
\end{proof}
\begin{crl}\label{crl: covering space theory}
If $d\equiv 0\pmod 2$ and $1\leq r\leq d-1$,
the inclusion map $S_r\hookrightarrow\pomod{d}$ induces an epimorphism
$$
\pi_1(S_r)\to \pi_1(\pomod{d})
\quad
\mbox{given by }\ \ \Z\oplus \Z \to \Z;
\ \ (a,b)\mapsto \pm b.
$$
\end{crl}
\begin{proof}
Note that the inclusion $S_r \stackrel{\subset}{\longrightarrow} \pomod{d}$ can be identified with the following composite of maps
\begin{equation}
\begin{CD}
S_r @>>\cong> \C^*\times \R^{d-r}\times \pomod{r} @>>> \pomod{d}
\\
@. (u,A,x) @>>> u\cdot\Phi_A(x)
\end{CD}
\end{equation}
Here $u\in\C^*$ acts by multiplication on the complex polynomial
formed by the last two coordinates. The loop in the $\C^*$ factor
is an orbit loop and maps trivially on $\pi_1$ by Lemma \ref{lmm: 3.9},
since $d$ is even. The $\R^{d-r}$ factor is contractible, while
Lemma \ref{lmm: pi1 induced} gives an isomorphism on the last factor.
This proves the assertion.
\end{proof}
\begin{definition}
{\rm
For each $d\geq 2$,
let $v_d:\widetilde{\pomod{d}} \to \pomod{d}$
denote the universal covering of $\pomod{d}$.
For each $0\leq r\leq d-1$, let $\widetilde{F}_r$ and $\widetilde{S}_r$ denote the subspaces of
$\widetilde{\pomod{d}}$ 
defined by
\begin{equation}
\widetilde{F}_r=v_d^{-1}(F_r)
\quad
\mbox{and}\quad
\widetilde{S}_r=v_d^{-1}(S_r).
\end{equation}
}
\end{definition}
\begin{remark}
{\rm 
Let $d\equiv 0$ (mod $2$).
\par
(i) Recall the homeomorphism
$
S_r\cong \C^*\times \R^{d-r}\times \pomod{r}
$
as in (\ref{eq: homeo-Sr}).
Then
it follows from 
Corollary \ref{crl: covering space theory} and covering-space classification that
there is a natural homeomorphism
\begin{equation}\label{eq: lifted degree stratum}
\widetilde S_r\cong
\C^*\times\R^{d-r}\times\widetilde{\pomod{r}}\qquad(r\geq1).
\end{equation}
The
 deck generator acts as a generator or its inverse
on the last factor and trivially on the other factors. 
This product decomposition concerns the pull-back covering;
 the $\C^*$ factor is not universally covered.
For $r=0$, $S_0\cong\C^*\times\R^d$ and its orbit loop maps
trivially on $\pi_1$ by Lemma \ref{lmm: 3.9}. Thus the restricted cover $\widetilde{S}_0$ is a disjoint union of copies of $S_0$ indexed
by $\Z$.
\par
In particular,
by using (\ref{eq: lifted degree stratum}), we easily obtain the  isomorphism
\begin{equation}\label{eq: S=po +po}
H_j(\widetilde{S}_r;\Z)\cong H_j(\widetilde{\pomod{r}};\Z)\oplus
H_{j-1}(\widetilde{\pomod{r}};\Z)
\end{equation}
for $1\leq r\leq d-1$ and $j\geq 0$, as an isomorphism of
$\Z[t,t^{-1}]$-modules. Negative-degree homology is understood to be zero.
In particular, $H_0(\widetilde S_r;\Z)\cong\Z$ has trivial deck action
for $r\geq1$.
\par\vspace{1mm}\par
(ii)
Let $x=(f,g,h)\in F_{r-1}$ and write $q=g+\sqrt{-1}h$.
\par
Since $F_{r}$ is open in the real linear coefficient space $L_r=\R^d\times \C^{r+1}$, there is a continuous positive radius function
$\delta (x)$ sufficiently small such that adding a coefficient $v$ of $z^r$ to $q$ remains in $F_r$ whenever
$|v|<\delta (x)$.\footnote{%
For example,  we can take $\delta (x) $ as a sufficiently small multiple of the distance to $L_r\setminus F_r$.
}
Thus,  there is an open tubular neighbourhood 
\begin{equation}
V_r\cong F_{r-1}\times D^2_{\rm open}
\end{equation}
of $F_{r-1}$ (in $F_r$) with zero section $F_{r-1}$.
This normal framing of $F_{r-1}\subset F_r$ pulls back to a deck-invariant framing
of  $\widetilde{F}_{r-1}\subset \widetilde{F}_r$.
Thus, there is a deck-equivariant homeomorphism
\begin{equation}
\widetilde{V}_r=v_d^{-1}(V_r)\cong \widetilde{F}_{r-1}\times D^2_{\rm open}
\end{equation}
with trivial action on the disk.
\par
Since $\widetilde{F}_r=\widetilde{V}_r\cup \widetilde{S}_r$, by excision and the relative cross product
there is an isomorphism
\begin{align*}
H_j(\widetilde{F}_r,\widetilde{S}_r;\Z)
&\cong
H_j(\widetilde{F}_{r-1}\times D^2_{\rm open},\widetilde{F}_{r-1}\times(D^2_{\rm open}\setminus\{0\});\Z)
\cong H_{j-2}(\widetilde{F}_{r-1};\Z).
\end{align*}
Thus, we obtain the following exact sequence of $\Z[t,t^{-1}]$-modules:
\begin{equation}\label{eq: main exact sequence 2}
\cdots \to H_j(\widetilde{S}_r;\Z) \to H_j(\widetilde{F}_r;\Z)
\to H_{j-2}(\widetilde{F}_{r-1};\Z)
\to \cdots 
\quad\qed
\end{equation}
}
\end{remark}
\paragraph{3.5 Proof of Theorem \ref{thm: simple up to dimension}.}
\par
Now we can give the proof of Theorem \ref{thm: simple up to dimension}.
\par\vspace{2mm}\par
\begin{proof}[Proof of Theorem \ref{thm: simple up to dimension}]
For odd $d$, the product decomposition (\ref{eq: univ. cover for d:odd})
gives trivial deck action in every degree, including $d=1$ by
(\ref{eq: conditions 31}). We prove (ii) by induction on even $d$.
Degree zero has trivial action because the universal cover is connected;
for $d=2$, simple connectedness gives $H_1(\widetilde{\pomod{2}};\Z)=0$.

Let $d\geq4$ be even and suppose the assertion holds for every smaller
even positive degree. Fix $1\leq k<d$ and set
\[
j_r=k-2(d-1-r)\qquad(0\leq r\leq d-1).
\]
Since $j_0<0$, $H_{j_0}(\widetilde F_0;\Z)=0$, and likewise for
$\widetilde S_0$. For $1\leq r\leq d-1$ use (\ref{eq: S=po +po}).
If $r$ is odd, both summands have trivial deck action by the odd-degree
product decomposition. If $r$ is even, then $r\leq d-2$ and
\[
j_r\leq 2r-d+1\leq r-1.
\]
Both summands have nilpotent action by the induction hypothesis
(with trivial action in degree zero and zero groups in negative degrees).
Replacing a deck generator by its inverse preserves nilpotence.
Consequently $H_{j_r}(\widetilde S_r;\Z)$ is nilpotent for every $r$.

Now induct on $r$ in the exact sequence (\ref{eq: main exact sequence 2}),
using $j_r-2=j_{r-1}$:
\[
H_{j_r}(\widetilde S_r;\Z)\longrightarrow
H_{j_r}(\widetilde F_r;\Z)\longrightarrow
H_{j_{r-1}}(\widetilde F_{r-1};\Z).
\]
The two outside modules are nilpotent, so the middle module is nilpotent
as well. At $r=d-1$ this proves nilpotence of
$H_k(\widetilde{\pomod{d}};\Z)$, as required.
\end{proof}
\begin{lemma}\label{lmm: nilpotent covering homology}
Let $N\geq 2$ and let
$f:X\to Y$ be a based continuous map between connected spaces $X$ and $Y$
such that $\pi_1(X)\cong \pi_1(Y)\cong \Z$.
Moreover, let $\widetilde{X}$ and $\widetilde{Y}$ denote the universal coverings of $X$ and $Y$, and
the following three conditions hold:
\begin{enumerate}
\item[$($a$)$]
The map $f$ is a homology equivalence up to dimension $N$.
\item[$($b$)$]
The deck action on $H_k(\widetilde{X};\Z)$ is nilpotent for any $2\leq k<N$.
\item[$($c$)$]
The deck action on $H_k(\widetilde{Y};\Z)$ is nilpotent for any $2\leq k\leq N$.
\end{enumerate}
Then the map $f$ is a homotopy equivalence up to dimension $N$.
\end{lemma}
\begin{remark}
{\rm 
Under the hypotheses of Lemma \ref{lmm: nilpotent covering homology},
the homomorphism 
\begin{equation}\label{eq: iso f pi1}
f_*:\pi_1(X)\stackrel{\cong}{\longrightarrow} \pi_1(Y)
\end{equation}
is an isomorphism.
Indeed, this easily follows from
the following commutative  diagram
\begin{equation}
\begin{CD}
\pi_1(X) @>f_*>> \pi_1(Y)
\\
@V{h}V{\cong}V  @V{h}V{\cong}V
\\
H_1(X;\Z) @>f_*>\cong> H_1(Y;\Z)
\end{CD}
\end{equation}
}
\end{remark}
\par
\begin{proof}[Proof of Lemma \ref{lmm: nilpotent covering homology}]
Use CW approximations and replace $f$ by its cellular mapping
cylinder inclusion. The lift of $X$ inside the universal cover of $M_f$ is
its universal cover because $f_*$ is a fundamental group isomorphism.
Put
\[
C=C_*(\widetilde M_f,\widetilde X;\Z),\qquad L_j=H_j(C),
\]
using cellular chains. The equivariant homology sequence gives the following exact sequence
\[
\begin{split}
0\longrightarrow&
\operatorname{coker}\bigl(H_j(\widetilde X)\to H_j(\widetilde Y)\bigr)
\longrightarrow L_j\\
\longrightarrow&
\ker\bigl(H_{j-1}(\widetilde X)\to H_{j-1}(\widetilde Y)\bigr)
\longrightarrow0.
\end{split}
\]
The hypotheses make both outside modules, and hence $L_j$,
nilpotent for $2\leq j\leq N$. Here $H_1$ of either cover is zero, and
$L_0=L_1=0$ because the covers are connected and simply connected.

The relative cellular chain groups are free $\Lambda$-modules, so
there is a short exact sequence of chain complexes
\[
0\longrightarrow C\xrightarrow{t-1}C
\longrightarrow C\otimes_\Lambda\Z\longrightarrow0.
\]
The last complex computes ordinary $H_*(M_f,X;\Z)$, which vanishes
through degree $N$ by the homology hypothesis. The resulting long
exact sequence shows that $t-1:L_j\to L_j$ is surjective for
$j\leq N$. Hence $L_j=(t-1)L_j$, and finite nilpotence implies
$L_j=(t-1)^aL_j=0$ for some integer $a\geq1$.
\par
Thus, we see that
$H_k(\widetilde M_f,\widetilde X;\Z)=0$
for any $k\leq N$.
Since $(\widetilde M_f,\widetilde X)$ is the
simply connected covering pair,
by the relative Hurewicz theorem \cite[\S 4.2]{Hatcher}, applied in the first possible nonzero degree,
we see that $(\widetilde M_f,\widetilde X)$ is $N$-connected.
The inclusion $\widetilde X\hookrightarrow\widetilde M_f$ therefore
induces isomorphisms on $\pi_k$ for $k<N$ and an epimorphism for
$k=N$, by the homotopy exact sequence of the pair. The covering
projections induce isomorphisms on higher homotopy groups, and
$f_*$ is an isomorphism on $\pi_1$ by (\ref{eq: iso f pi1}).
Since the mapping cylinder retracts onto $Y$, this proves the assertion.
\end{proof}

\begin{crl}\label{crl: simple even case}
If $d\equiv0$ $({\rm mod}\ 2)$, the space $\po^{d,3}_1(\R)$ is
simple up to dimension $d$.
\end{crl}

\begin{proof}
For $d=2$, the assertion says only that the fundamental group acts
trivially on itself; this holds since it is isomorphic to $\Z$.
Suppose now that $d\geq4$ is even, and consider the stabilization map
\[
s^{d-1,3}_{1,\R}:\po^{d-1,3}_1(\R)\longrightarrow
\po^{d,3}_1(\R).
\]
This map induces a fundamental group isomorphism and an ordinary
homology equivalence through dimension $d-1$, by Lemma
\ref{lmm: KY13 (m,n)=(3,1)}. The odd-degree source space $\po^{d-1,3}_1(\R)$ has trivial deck
action on the homology of its universal cover in every degree, by
Lemma \ref{lmm: polyd: odd}. The target space $\po^{d,3}_1(\R)$ has nilpotent deck action
through degree $d-1$, by Theorem \ref{thm: simple up to dimension}.
By  Lemma \ref{lmm: nilpotent covering homology} with $N=d-1$, we see that
$s^{d-1,3}_{1,\R}$ is a homotopy equivalence up to dimension $d-1$.
\par
Let $1\leq k<d$ and let
$(\alpha,\beta)\in \pi_1(\po^{d,3}_1(\R))\times \pi_k(\po^{d,3}_1(\R))$ be any element.
Since $s^{d-1,3}_{1,\R}$ is a homotopy equivalence up to dimension $d-1$, there exists a pair
$(a,b)\in \pi_1(\po^{d-1,3}_1(\R))\times \pi_k(\po^{d-1,3}_1(\R))$ such that
$((s^{d-1,3}_{1,\R})_*(a),(s^{d-1,3}_{1,\R})_*(b))=(\alpha,\beta)$.
Since $d-1\equiv 1$ (mod $2$), $\po^{d-1,3}_1(\R)$ is simple by Corollary \ref{crl: Polyd,3,1 simple}.
So $a\cdot b=b$.
Thus,
\[
\alpha\cdot\beta
=(s^{d-1,3}_{1,\R})_*(a)\cdot (s^{d-1,3}_{1,\R})_*(b)=(s^{d-1,3}_{1,\R})_*(a\cdot b)=
(s^{d-1,3}_{1,\R})_*(b)=\beta.
\]
Thus, $\po^{d,3}_1(\R)$ is simple up to dimension $d$.
\end{proof}

\section{Proofs of the main results}\label{section: proofs}
\paragraph{4.1 Proofs of the main results.}
We first prove Theorems \ref{thm: KY14: homotopy simple} and
\ref{thm: KY14: (m,n)=(3,1)}.

\begin{proof}[Proof of Theorem \ref{thm: KY14: homotopy simple}]
For odd $d$, the assertion is Corollary
\ref{crl: Polyd,3,1 simple}, including $d=1$ by
$\po^{1,3}_1(\R)\simeq S^1$. For even $d$, it is Corollary
\ref{crl: simple even case}.
\end{proof}

\begin{proof}[Proof of Theorem \ref{thm: KY14: (m,n)=(3,1)}]
For $d=1$, both maps induce isomorphisms on fundamental groups by
Lemma \ref{lmm: KY13 (m,n)=(3,1)}, which is more than the required
epimorphism in the endpoint dimension. All spaces are path-connected.
Assume henceforth that $d\geq2$.

We apply Lemma \ref{lmm: nilpotent covering homology} directly,
without using Corollary \ref{crl: simple even case}. In either
parity, the deck action on the source universal-cover homology is
nilpotent for $k<d$, by Theorem
\ref{thm: simple up to dimension}. For the natural map, the target
has the homotopy decomposition
\[
(\Omega^2_d\CP^2)^{\Z_2}
\simeq \Omega^2S^5\times\Omega S^3\times S^1.
\]
The first two factors are simply connected. Thus its universal
cover has trivial deck action on homology in every degree.

For the stabilization, if $d$ is even, its target has odd degree
and the deck action on its universal-cover homology is trivial in
every degree. If $d$ is odd, its target has even degree $d+1$ and
Theorem \ref{thm: simple up to dimension} gives nilpotent deck action
in every degree $k<d+1$, hence in particular for all $k\leq d$.

Both maps induce fundamental group isomorphisms by Lemma
\ref{lmm: KY13 (m,n)=(3,1)} and ordinary homology isomorphisms
through dimension $d$ by Theorem \ref{thm: KY13}. All hypotheses
of Lemma \ref{lmm: nilpotent covering homology}, with $N=d$, are
therefore satisfied. It follows that both maps induce isomorphisms
on $\pi_k$ for $k<d$ and epimorphisms on $\pi_d$, as required.
No hypothesis on the action in source degree $d$, and no
injectivity conclusion in that degree, is asserted or used.
\end{proof}

\paragraph{4.2 The space $(\Omega^2\CP^N)^{\Z_2}$.}
Finally, we give the proof of the homotopy equivalence (\ref{eq: 1.7}).
For this purpose, we recall the following facts.
\par\vspace{1.5mm}\par
From now on,
we regard $S^2$ as the union $S^2=D^2_+\cup D^2_{-}$, 
where  
let
$D^2_+$ and $D^2_{-}$ denote the northern and southern hemispheres given by
\begin{equation}\label{eq: S2 decomposition}
\begin{cases}
D^2_+&=\{(x,y,z)\in \R^3:x^2+y^2+z^2=1,\ z\geq 0\},
\\
D^2_{-}&=\{(x,y,z)\in \R^3:x^2+y^2+z^2=1,\ z\leq 0\}.
\end{cases}
\end{equation}
We identify the intersection $D^2_+\bigcap D^2_{-} = \{(x,y,0)\in \R^3: x^2+y^2=1\}$ with the unit circle $S^1$.
Note that the group $\Z_2=\{\pm 1\}$ acts on $S^2$ by reflection in the equator
(i.e. 
$S^2\ni x=(x_1,y_1,z_1)  \mapsto \overline{x}=(x_1,y_1,-z_1)$) and on $\CP^N$
by complex conjugation.
\par
\begin{remark}\label{remark: space of equiv. maps}
{\rm 
Let $N\geq 2$.
\par
(i)
Consider the restriction map
\begin{equation}
r_{D^2_+}:(\Omega^2\CP^N)^{\Z_2}\to \Map^*(D^2_+,S^1;\CP^N,\RP^N)
\end{equation}
given by $r_{D^2_+}(f)=f\vert D^2_+$
for $f\in (\Omega^2\CP^N)^{\Z_2}$.
\par
Restriction to $D^2_+$ identifies the fixed point space
$(\Omega^2\CP^N)^{\Z_2}$ with the space
$\Map^*(D^2_+,S^1;\CP^N,\RP^N)$ of based relative maps
from $(D^2_+,S^1)$ to $(\CP^N,\RP^N)$.
More precisely, consider the  map
\begin{equation}
\mathcal{D}:\Map^*(D^2_+,S^1;\CP^N,\RP^N)\to (\Omega^2\CP^N)^{\Z_2}
\end{equation}
defined by\footnote{%
The condition that \(u(S^1)\subset \RP^N\) is precisely what makes this doubled
map $\mathcal{D}(u)$
continuous along the equator.
}
\begin{equation*}
\mathcal{D}(u)(x)=
\begin{cases}
u(x)  & \mbox{ if }x\in D^2_+
\\
\overline{u(\overline{x})} & \mbox{ if }x\in D^2_{-}
\end{cases}
\ \ 
\mbox{for $(u,x)\in \Map^*(D^2_+,S^1;\CP^N,\RP^N)\times S^2$.}
\end{equation*}
It is easy to see that
$r_{D^2_+}\circ \mathcal{D}=\mbox{id}$ and
$\mathcal{D}\circ r_{D^2_+}=\mbox{id}$.
Thus, this doubling map
\begin{equation}\label{eq: equiv. double map}
\mathcal{D}:\Map^*(D^2_+,S^1;\CP^N,\RP^N)\stackrel{\cong}{\longrightarrow} (\Omega^2\CP^N)^{\Z_2}
\end{equation}
is a homeomorphism whose inverse is $r_{D^2_+}$.
\par
Consequently, the path components of the fixed point space
\((\Omega^2_d\CP^N)^{\Z_2}\) are identified with those relative homotopy classes
in
\begin{equation}\label{eq: identification}
     \pi_0(\Map^*(D^2_+,S^1;\CP^N,\RP^N))  
     \cong \pi_2(\CP^N,\RP^N)
\end{equation}
whose doubled sphere has degree $d$.\footnote{%
The obvious relative homeomorphism $(D^2_+,S^1)\cong (D^2,S^1)$
gives rise to a homeomorphism
$\Map^*(D^2_+,S^1;\CP^N,\RP^N)\cong \Map^*(D^2,S^1;\CP^N,\RP^N)$.
By  means of this identification, we obtain the isomorphism
(\ref{eq: identification}).
}
\par
(ii)
We now describe these relative classes using the Maslov index.
\par
Since
\(\RP^N\subset \CP^N\) is the fixed point set of complex conjugation, which
is anti-symplectic for the Fubini--Study form, it is Lagrangian, and in
particular totally real. Hence the Maslov class \cite[page 5]{KS}
\[
        \mu_{\RP^N}\in H^2(\CP^N,\RP^N;\Z)
\]
and the corresponding Maslov index of a relative disc are defined. As
recorded by Konstantinov and Smith \cite[Section 3.6]{KS}, the group
\begin{equation}
        \pi_0(\Map^*(D^2_+,S^1;\CP^N,\RP^N))
        \cong \pi_2(\CP^N,\RP^N)\cong \Z
\end{equation}
is generated by the class \(\beta\) of the standard upper hemisphere in
\[
        (\CP^1,\RP^1)\subset (\CP^N,\RP^N),
\]
and this generator has Maslov index \(N+1\). Thus every relative class is
of the form \(k\beta\), and its Maslov index is \(k(N+1)\).
\par
Moreover,
the doubled sphere \(\mathcal D_*(\beta)\) is represented by the standard line
$
        \CP^1\hookrightarrow \CP^N,
$
and therefore has degree \(1\). 
Thus, if $\iota :(\Omega^2\CP^N)^{\Z_2}\hookrightarrow \Omega^2\CP^N$ denotes the inclusion,
the composite
\begin{equation*}
\begin{CD}
\pi_0(\Map^*(D^2_+,S^1;\CP^N,\RP^N)) @>\mathcal{D}_*>\cong> \pi_0((\Omega^2\CP^N)^{\Z_2})
@>\iota_*>> \pi_0(\Omega^2\CP^N)
\\
@AA{\cong}A @. @V{\deg}V{\cong}V
\\
\pi_2(\CP^N,\RP^N)= \Z \cdot \beta
@. @. \Z
\end{CD}
\end{equation*}
sends \(\beta\) to \(1\). 
Since $\pi_2(\CP^N,\RP^N)=\Z\cdot \beta$,
the composite
$\deg\circ \iota_*\circ \mathcal{D}_*$ 
is 
an isomorphism. 
Furthermore, since $\mathcal{D}_*$ and
$\deg:\pi_0(\Omega^2\CP^N)\stackrel{\cong}{\longrightarrow}\Z$ are isomorphisms, it follows that
\begin{equation}
\iota_*:\pi_0((\Omega^2\CP^N)^{\Z_2})
\stackrel{\cong}{\longrightarrow}
\pi_0(\Omega^2\CP^N)
\end{equation}
is an isomorphism.
Consequently, 
for each
\(d\in \Z\), there is exactly one relative homotopy class whose double has
degree \(d\). Therefore, the space
$
        (\Omega^2_d\CP^N)^{\Z_2}
$
is path-connected for each $d\in \Z$.
\par
Furthermore, we also obtain the following path-component decomposition:
\begin{equation}
(\Omega^2\CP^N)^{\Z_2}=\bigcup_{d\in \Z}(\Omega^2_d\CP^N)^{\Z_2}
\quad
\mbox{(disjoint union)}.
\qquad
\qed
\end{equation}
}
\end{remark}
\begin{definition}
{\rm
Let $N\geq 2$.
For each $d\in \Z$, let $M_d$ denote the subspace 
of $\Map^*(D^2_+,S^1;\CP^N,\RP^N)$
defined by
\begin{equation}\label{def: Md}
M_d=\{u\in \Map^*(D^2_+,S^1;\CP^N,\RP^N):\mathcal{D}(u)\in (\Omega^2_d\CP^N)^{\Z_2}\}.
\end{equation}
}
\end{definition}
\begin{remark}
{\rm
(i)
Since $\mathcal{D}$ is a homeomorphism by (\ref{eq: equiv. double map})
and $(\Omega^2_d\CP^N)^{\Z_2}$ is path-connected,
the space $M_d$ is path-connected. 
Moreover, we obtain the following decomposition of path components:
\begin{equation}
 \Map^*(D^2_+,S^1;\CP^N,\RP^N)=\bigcup_{d\in \Z}M_d
 \quad
\mbox{(disjoint union)}.
\end{equation}
\par
(ii)
Under the identification
$\pi_0(\Map^*(D^2_+,S^1;\CP^N,\RP^N))\cong\Z$,
the path component $M_d$ corresponds to $d$. Thus, we have the correspondence
\begin{equation}\label{eq: correspondence}
\pi_0(\Map^*(D^2_+,S^1;\CP^N,\RP^N)) 
\cong
\Z ,\quad M_d\leftrightarrow d.
\end{equation}
We also remark that the restriction
\begin{equation}\label{eq: mathcal Dd}
\mathcal{D}\vert M_d:M_d
\stackrel{\cong}{\longrightarrow}
(\Omega^2_d\CP^N)^{\Z_2}
\end{equation}
is a homeomorphism for each $d\in \Z$.
\qed
}
\end{remark}
\begin{definition}
{\rm
Let $N\geq 2$.
For each $\epsilon\in \{0,1\}$, let
\begin{equation}
\Map^*_{\epsilon}=\Map^*_{\epsilon}(D^2_+,S^1;\CP^N,\RP^N)
\end{equation}
denote the subspace of
$\Map^*(D^2_+,S^1;\CP^N,\RP^N)$
given by
\begin{align}
\Map^*_{\epsilon}&=\Map^*_{\epsilon}(D^2_+,S^1;\CP^N,\RP^N)
\\
\nonumber
&=
\{f\in \Map^*(D^2_+,S^1;\CP^N,\RP^N):
f\vert S^1\in \Omega_{\epsilon}\RP^N\}.
\end{align}
}
\end{definition}
\begin{lemma}\label{lmm: Md}
If $N\geq 2$,
the two subspaces have the following decompositions into path components:
\begin{equation}\label{eq: path-component even or odd}
\begin{cases}
\dis \Map^*_{0}(D^2_+,S^1;\CP^N,\RP^N)=\bigcup_{k\in \Z}M_{2k}
& (\mbox{disjoint union}),
\\
\dis \Map^*_{1}(D^2_+,S^1;\CP^N,\RP^N)=\bigcup_{k\in \Z}M_{2k+1}
& (\mbox{disjoint union}).
\end{cases}
\end{equation}
\end{lemma}
\begin{proof}
Let
$r_{S^1}:\Map^*(D^2_+,S^1;\CP^N,\RP^N)\to \Omega \RP^N$
denote the restriction map given by
$r_{S^1}(u)=u\vert S^1$ for
$u\in \Map^*(D^2_+,S^1;\CP^N,\RP^N)$, and
let $q:(D_+^2,S^1)\to (S^2,*)$ be the pinch map.
\par
Then
by  \cite[Lemma 3.2]{KY13} for 
$(X,A;Y,B)=(D^2_+,S^1;\CP^N,\RP^N)$, we obtain the following fibration sequence
\begin{equation}\label{eq: restrict fib. seq}
\Omega^2\CP^N 
\stackrel{q^{\#}}{\longrightarrow}
\Map^*(D^2_+,S^1;\CP^N,\RP^N)
\stackrel{r_{S^1}}{\longrightarrow}
\Omega \RP^N,
\end{equation}
where
the map $q^{\#}$ is given by $q^{\#}(f)=f\circ q$ for $f\in \Omega^2\CP^N$.
\par
Consider the exact sequence of the pair $(\CP^N,\RP^N)$. Its boundary operator
$\partial_2:\pi_2(\CP^N,\RP^N)\to \pi_1(\RP^N)$ is given by
$\partial_2([f])=[f\vert S^1]$.
Thus, if we
consider the exact homotopy sequence of the fibration
(\ref{eq: restrict fib. seq}), then 
under the standard identifications
\[
\begin{split}
\pi_0(\Omega^2\CP^N)&=\pi_2(\CP^N)\cong\Z,\\
\pi_0(\Map^*)&=\pi_2(\CP^N,\RP^N)\cong\Z,\\
\pi_0(\Omega\RP^N)&=\pi_1(\RP^N)\cong\Z/2,
\end{split}
\]
we have the following commutative  diagram of the exact sequence.
$$
\begin{CD}
\pi_2(\CP^N,\RP^N) @>\partial_2>> \pi_1(\RP^N) @>>> 0=\pi_1(\CP^N)
\\
@VV{\cong}V @VV{\cong}V @.
\\
\pi_0(\Map^*) @>(r_{S^1})_*>> \pi_0(\Omega \RP^N) @. 
\end{CD}
$$
where we set $\Map^*=\Map^*(D^2_+,S^1;\CP^N,\RP^N)$.
\par
It is easy to see that  $(r_{S^1})_*$ is an epimorphism. 
More precisely,
with the above identification
this epimorphism 
$(r_{S^1})_*:\Z\to \Z/2$ is reduction mod $2$.
Therefore,
the relevant part of this sequence (induced from (\ref{eq: restrict fib. seq}))
is the following exact sequence:
\begin{equation*}
\begin{CD}
\pi_0(\Omega^2\CP^N)
  @>{(q^{\#})_*}>>
  \pi_0(\Map^*)
  @>{(r_{S^1})_*}>>
  \pi_0(\Omega\RP^N) @>>> 0
  \\
@V{\cong}VV @V{\cong}VV @V{\cong}VV @.\\
\Z @>{}>> \Z @>>> \Z/2 @.
\end{CD}
\end{equation*}
Under the above identifications,
$\mbox{Ker }(r_{S^1})_*=\mbox{Im }(q^{\#})_*=2\Z$.
Thus, $(q^{\#})_*$ is multiplication by $\pm 2$; after choosing the
generator of $\pi_0(\Map^*)\cong \Z$ compatibly, it is multiplication by
$2$.
Hence, we obtain the following commutative diagram with exact rows.
{\small
\begin{equation}\label{CD: x 2 diagram}
\begin{CD}
0 @>>> \pi_0(\Omega^2\CP^N)
  @>{(q^{\#})_*}>>
  \pi_0(\Map^*)
  @>{(r_{S^1})_*}>>
  \pi_0(\Omega\RP^N) @>>> 0\\
@. @V{\cong}VV @V{\cong}VV @V{\cong}VV @.\\
@. \Z @>{\times 2}>> \Z @>>> \Z/2 @.
\end{CD}
\end{equation}
}
\par\noindent
For $u\in\Map^*$, let $[u]\in\pi_0(\Map^*)$ denote its path component.
Then
$$
\begin{cases}
(r_{S^1})_*([u])=0
\Leftrightarrow u\vert S^1\in \Omega_0\RP^N 
\Leftrightarrow u\in 
\Map^*_0(D^2_+,S^1;\CP^N,\RP^N),
\\
(r_{S^1})_*([u])\not=0
\Leftrightarrow 
u\vert S^1\in \Omega_1\RP^N
\Leftrightarrow u\in 
\Map^*_1(D^2_+,S^1;\CP^N,\RP^N).
\end{cases}
$$
Under the identification
$\pi_0(\Map^*(D^2_+,S^1;\CP^N,\RP^N))\cong\Z$,
the diagram (\ref{CD: x 2 diagram}) shows that
$$
\begin{cases}
\pi_0(\Map^*_0(D^2_+,S^1;\CP^N,\RP^N))=\Z_{\rm even}=\{2k:k\in \Z\},
\\
\pi_0(\Map^*_1(D^2_+,S^1;\CP^N,\RP^N))=\Z_{\rm odd}=\{2k+1:k\in \Z\}.
\end{cases}
$$
The correspondence (\ref{eq: correspondence}) now gives the desired
decompositions into path components (\ref{eq: path-component even or odd}).
\end{proof}
\begin{lemma}\label{lmm: Lemma 3.3 in KY13}
If $N\geq 2$, there is a homotopy equivalence
\begin{equation}\label{eq: homotopy equiv d=0}
(\Omega_d^2\CP^N)^{\Z_2}\simeq \Omega^2S^{2N+1}\times \Omega S^N
\quad
\mbox{for every }\ d\in \Z.
\end{equation}
\end{lemma}
\begin{proof}
Choose the loop multiplication $*$ to be induced by a
$\Z_2$-equivariant pinch map $S^2\to S^2\vee S^2$. Then the product
of two equivariant maps is again equivariant.
By Remark \ref{remark: space of equiv. maps} (ii), the space
$(\Omega^2_d\CP^N)^{\Z_2}$ is nonempty. Choose 
$g\in(\Omega^2_d\CP^N)^{\Z_2}$ and define the map
\[
A_g:(\Omega^2_0\CP^N)^{\Z_2}\longrightarrow
(\Omega^2_d\CP^N)^{\Z_2}
\ 
\mbox{ by }
\  A_g(f)=f*g \ \mbox{ for }\ f\in (\Omega^2_0\CP^N)^{\Z_2}.
\]
\par
We claim that this map is a homotopy equivalence
\begin{equation}\label{eq: d=0}
A_g:
(\Omega^2_0\CP^N)^{\Z_2}\stackrel{\simeq}{\longrightarrow} (\Omega^2_d\CP^N)^{\Z_2}.
\end{equation}
Indeed, 
choose 
$
h\in(\Omega^2_{-d}\CP^N)^{\Z_2}
$
and
define the map
$$
A_h:(\Omega^2_d\CP^N)^{\Z_2}\longrightarrow
(\Omega^2_0\CP^N)^{\Z_2}
\ 
\mbox{ by $A_h(f)=f*h$ \ for \ $f\in (\Omega^2_d\CP^N)^{\Z_2}$.}
$$
Then,
$
g*h,\ h*g\in(\Omega^2_0\CP^N)^{\Z_2}.
$
Since this space $(\Omega^2_0\CP^N)^{\Z_2}$ is path-connected, both $g*h$ and $h*g$ can be joined
to the constant map by paths of equivariant maps of degree zero. Using the
homotopy associativity and homotopy unit of the loop product, these paths
induce homotopies
\[
A_h\circ A_g \simeq \mbox{\rm id}_{(\Omega^2_0\CP^N)^{\Z_2}},
\qquad
A_g\circ A_h \simeq \mbox{\rm id}_{(\Omega^2_d\CP^N)^{\Z_2}}.
\]
Hence $A_g$ is a homotopy equivalence with homotopy inverse $A_h$.
%
\par
It therefore suffices to prove (\ref{eq: homotopy equiv d=0}) for $d=0$.
\par
Let  $j_N:\RP^N\to \CP^N$ 
and $i_N:S^N\to S^{2N+1}$
denote the natural inclusions.
Let $\xi_N$ denote the restriction fibration
\[
\rho_N
:\Map^*(D^2,\CP^N)\longrightarrow \Omega \CP^N
\]
given by
$\rho_N(f)=f\vert S^1$.
Since $\Map^*(D^2,\CP^N)$ is contractible, the fiber of 
$\rho_N$ 
is
$\Omega^2\CP^N$. 
We denote by $\Map^{**}_0=\Map^*_0(D^2,S^1;\CP^N,\RP^N)$ the space given by
\begin{align*}
\Map^{**}_0&=
\Map^*_0(D^2,S^1;\CP^N,\RP^N)
\\
&=\{f\in \Map^*(D^2,S^1;\CP^N,\RP^N):
f\vert S^1\in \Omega_0 \RP^N\}.
\end{align*}
\par
Now consider the pullback fibration  $(\Omega j_N)^* \xi_N$ of $\xi_N$ by the map
$\Omega j_N: \Omega_0\RP^N\to \Omega \CP^N$.
If $E_N$ denotes the total space of the fibration $(\Omega j_N)^* \xi_N$, 
we obtain the following commutative diagram of fibration sequences:
$$
\begin{CD}
\Omega^2\CP^N @>>> E_N @>>> \Omega_0\RP^N
\\
\Vert @. @VVV  @V{\Omega j_N}VV
\\
\Omega^2\CP^N @>>> \Map^*(D^2,\CP^N) @>\rho_N>> \Omega \CP^N
\end{CD}
$$
%
Hence, the space $E_N$ can be expressed as follows:
\begin{align*}
E_N&=
\{(f,\omega)\in \Map^*(D^2,\CP^N)\times \Omega_0\RP^N:
\rho_N(f)=(\Omega j_N)(\omega)\}
\\
&= \{(f,\omega)\in \Map^*(D^2,\CP^N)\times \Omega_0\RP^N : f(S^1)\subset \RP^N,
f\vert S^1=\omega \}
\\
&\cong 
\{f\in \Map^*(D^2,S^1;\CP^N,\RP^N):f\vert S^1\in \Omega_0\RP^N\}
\\
&=
\Map^*_0(D^2,S^1;\CP^N,\RP^N).
\end{align*}
Thus we obtain a fibration sequence
\begin{equation}\label{eq: fibration4}
\Omega^2\CP^N\to 
\Map^*_0(D^2,S^1;\CP^N,\RP^N)
\stackrel{r_{S^1}}{\longrightarrow}
\Omega_{0}\RP^N,
\end{equation}
where the map $r_{S^1}$ is given by $r_{S^1}(f)=f\vert S^1$ as in the proof of Lemma
\ref{lmm: Md}.
\par
Let
$\gamma_N:S^{2N+1}\to \CP^N$ and
$\gamma_{N,\R}:S^N\to \RP^N$ denote the Hopf fibration
and the usual double covering, respectively.
Consider the following commutative diagram
$$
\begin{CD}
\Omega S^N @>\Omega i_N>> \Omega S^{2N+1}
\\
@V{\Omega \gamma_{N,\R}}V{\simeq}V @V{\Omega\gamma_N}VV
\\
\Omega_{0} \RP^N @>\Omega j_N>> \Omega \CP^N
\end{CD}
$$
The inclusion $i_N:S^N\to S^{2N+1}$ is null-homotopic, whereas
$\Omega\gamma_{N,\R}:\Omega S^N\stackrel{\simeq}{\longrightarrow}\Omega_0\RP^N$ is a homotopy
equivalence. It follows that
$\Omega j_N:\Omega_0\RP^N\to\Omega\CP^N$ is null-homotopic.
Pullbacks of a fibration along homotopic maps are fiber-homotopy equivalent.
Consequently, $(\Omega j_N)^*\xi_N$ is fiber-homotopy equivalent to the
pullback of $\xi_N$ along the constant map. The latter is the product
fibration
\[
\Omega^2\CP^N\longrightarrow
\Omega^2\CP^N\times\Omega_0\RP^N
\longrightarrow\Omega_0\RP^N.
\]
Therefore, there is a homotopy equivalence
\[
\Map^*_0(D^2,S^1;\CP^N,\RP^N)
\simeq
\Omega^2\CP^N\times\Omega_0\RP^N.
\]
Since there is a homotopy equivalence 
$\Omega^2\CP^N\simeq \Omega^2S^{2N+1}\times \Z$,
we obtain a homotopy equivalence
\begin{align*}
\Map^*_0(D^2,S^1;\CP^N,\RP^N)
&\simeq \Omega^2\CP^N\times \Omega_0\RP^N
\simeq 
\Omega^2S^{2N+1}\times \Omega S^N\times \Z.
\end{align*}
By identifying $D^2$ with $D^2_+$, we also obtain the homotopy equivalence
\begin{equation}
\Map^*_0(D^2_+,S^1;\CP^N,\RP^N)\simeq \Omega^2S^{2N+1}\times \Omega S^N\times \Z.
\end{equation}
Since $\Omega^2S^{2N+1}\times\Omega S^N$ is path-connected, each path
component of the right-hand side has this homotopy type. It follows from
(\ref{eq: path-component even or odd}) that
\begin{equation}
M_{2k}\simeq \Omega^2S^{2N+1}\times \Omega S^N
\quad
\mbox{ for every }k\in \Z.
\end{equation}
By (\ref{eq: mathcal Dd}) we finally obtain the homotopy equivalence
\begin{equation}\label{eq: space of maps}
(\Omega^2_0\CP^N)^{\Z_2}\cong 
M_0\simeq \Omega^2S^{2N+1}\times \Omega S^N.
\end{equation}
This completes the proof.
\end{proof}
Using (\ref{eq: mathcal Dd}) and (\ref{eq: homotopy equiv d=0}), we also obtain the following result.
\begin{crl}
If $N\geq 2$, there is a homotopy equivalence
$$
M_d\simeq \Omega^2 S^{2N+1}\times \Omega S^N
\quad
\mbox{ for any }\ d\in\Z.
\qquad\qquad\qquad
\qed
$$
\end{crl}
\section{Low-dimensional computations}\label{section: optimality}

We conclude by determining the homotopy types in degrees $2$ and $3$
and showing that the stability range of Theorem \ref{thm: KY14: (m,n)=(3,1)}
is optimal in both degrees.

\begin{lemma}\label{lmm: optimality d=2}
There is a homotopy equivalence
\begin{equation}\label{eq: po231}
\po^{2,3}_1(\R)\simeq S^1\vee S^2\vee S^3.
\end{equation}
\end{lemma}
\begin{proof}
Write a triple of monic quadratic polynomials as a vector-valued polynomial
$F(z)=z^2 e+za+b$, where $e=(1,1,1)$ and $a,b\in\R^3$.
Put $P=e^{\perp}=\{x\in \R^3:\langle x,e\rangle=0\}
\cong\R^2$ and write
$a=\alpha e+B$, $b=\beta e+C$, with $B,C\in P$.
Substitution of $z=t-\alpha/2$ gives
\[
F(t-\alpha/2)=(t^2+D)e+tB+E,
\qquad D=\beta-\alpha^2/4,\quad E=C-(\alpha/2)B.
\]
This change of coefficient coordinates gives a homeomorphism
$\po^{2,3}_1(\R)\cong\R\times Y$, where $Y$ is the open subset
of $\R\times P\times P$ consisting of triples $(D,B,E)$ for which
\[
t^2+D=0,\qquad tB+E=0
\]
have no simultaneous solution $t\in\C$. If $B\ne0$, such a solution
would have to be real. 
If $B=0$, the common-root condition excludes
precisely $E=0$, independently of $D$.
\par
Consider the open cover $Y=U\cup V$ given by
\[
\begin{aligned}
U&=\{(D,B,E)\in Y:B\ne0\},\\
V&=\{(D,B,E):E\ne0,\ |B|^2(1+|D|)<|E|^2\}.
\end{aligned}
\]
The inequality implies that $V\subset Y$, since a common root with
$B\ne0$ would imply $|E|^2=-D|B|^2$. Every point of $Y$ with $B=0$
belongs to $V$.

Fix an orientation of $P$ and let
$J:P\to P$ be rotation through $90^\circ$. On $U$ write uniquely
$E=\lambda B+\mu JB$, and set $\delta=D+\lambda^2$.
The excluded locus is $(\mu,\delta)=(0,0)$, so
\[
U\cong(P\setminus\{0\})\times\R\times
(\R^2\setminus\{0\})\simeq S^1_B\times S^1_M.
\]
Scaling $B$ to zero, then $D$ to zero, and normalizing $E$ gives
$V\simeq S^1_E$. 
\par
On $W=U\cap V$, the directions of $B$ and $E$
give $W\simeq S^1_B\times S^1_E$, since the remaining coordinates
$D$, $|E|$ and $|B|/|E|\in(0,(1+|D|)^{-1/2})$ form a contractible factor.

Under these identifications, $W\to V$ is represented by projection
onto $S^1_E$. If $\mu=0$ on $W$, then
$\lambda^2>1+|D|$ and hence $\delta>0$. Therefore the point
$(\mu,\delta)$ lies in the plane with the nonpositive $\delta$-axis
removed, a contractible space. The map $W\to U$ is thus represented
by projection onto $S^1_B$ followed by inclusion in the first factor.
Since $U,V$ form an open cover, $Y$ is homotopy equivalent to the
homotopy pushout of
\[
S^1_B\times S^1_M\longleftarrow S^1_B\times S^1_E
\longrightarrow S^1_E.
\]
Write $A=S^1_B$, $M=S^1_M$ and $K=S^1_E$. The homotopy pushout
of the two projections $A\leftarrow A\times K\to K$ is the join
$A*K\cong S^3$. By pasting homotopy pushouts, $Y$ is homotopy
equivalent to the homotopy pushout of
\[
A\times M\longleftarrow A\longrightarrow A*K.
\]
The map $A\to A*K$ is null-homotopic since $A$ is coned off in
the join. Thus
\[
Y\simeq 
\Big(\frac{S^1\times S^1}{S^1\times\{*\}}\Big)\vee S^3
\simeq S^1\vee S^2\vee S^3.
\]
The last equivalence follows from the usual CW structure on the torus,
since the commutator attaching map becomes null-homotopic after
collapsing one circle. Since $\R\times Y\simeq Y$, this proves the homotopy equivalence (\ref{eq: po231}).
\end{proof}

\begin{lemma}\label{lmm: the case d=3}
There are homotopy equivalences
\begin{equation}\label{eq: po331}
\po^{3,3}_1(\R)\simeq S^1\times(S^2\vee S^3),
\qquad
\widetilde{\po}^{3,3}_1\simeq S^2\vee S^3.
\end{equation}
\end{lemma}
\begin{proof}
Let $U=\widetilde{\po}^{3,3}_1$. By (\ref{eq: univ. cover for d:odd}),
$\po^{3,3}_1(\R)\simeq S^1\times U$. It follows from
(\ref{eq: stable homotopy equiv}) and (\ref{eq: D1, D2}) that
\[
\po^{3,3}_1(\R)\simeq_s S^1\vee S^2\vee S^3\vee S^3\vee S^4.
\]
The K\"unneth formula for $S^1\times U$ now gives
\[
H_k(U;\Z)\cong
\begin{cases}
\Z&\text{if }k=0,2,3,\\
0&\text{otherwise}.
\end{cases}
\]
Since $U$ is simply connected, the Hurewicz homomorphism
$h_2:\pi_2(U)\stackrel{\cong}{\longrightarrow} H_2(U;\Z)\cong \Z$ is an isomorphism. Choose a based map
$\varphi:S^2\to U$ representing a generator of $\pi_2(U)\cong\Z$.
Replacing it by its
mapping-cylinder inclusion, we may regard $(U,S^2)$ as a
$2$-connected pair. Relative Hurewicz gives an isomorphism
$\pi_3(U,S^2)\stackrel{\cong}{\longrightarrow} H_3(U,S^2;\Z)$. The homotopy and homology exact
sequences of this pair give the commutative diagram
\[
\begin{CD}
\pi_3(U) @>{\rm onto}>> \pi_3(U,S^2) @>>> \pi_2(S^2) @>{\varphi_*}>\cong> \pi_2(U)\\
@V{h_3}VV @VV{\cong}V @VV{\cong}V @V{h_2}V{\cong}V\\
H_3(U;\Z) @>>\cong> H_3(U,S^2;\Z) @>>> H_2(S^2;\Z) @>{\varphi_*}>\cong> H_2(U;\Z).
\end{CD}
\]
The first map in the top row is surjective because the last map is
an isomorphism. The first map in the bottom row is an isomorphism
because $H_3(S^2;\Z)=0$ and the last map is an isomorphism.
Hence $h_3$ is surjective. Choose a based map $\psi:S^3\to U$
whose Hurewicz image is a generator of $H_3(U;\Z)\cong \Z$.
The map $S^2\vee S^3\to U$ with restrictions $\varphi$ and $\psi$
induces an isomorphism on integral homology. Both spaces are simply
connected and of CW type, so the homological Whitehead theorem proves
the desired homotopy equivalence (\ref{eq: po331}).
\end{proof}
\begin{remark}\label{remark: best possible}
{\rm
(i)
Put $\Lambda=\Z[t,t^{-1}]$. Lemma \ref{lmm: optimality d=2} identifies
$\pi_2(\po^{2,3}_1(\R))$ with the regular $\Lambda$-module: the
universal cover of the wedge has one lift of $S^2\vee S^3$ at
each integer, and Hurewicz identifies its second homology with
its second homotopy group. In particular, the fundamental-group
action is nontrivial and the even-degree simplicity range of
Theorem \ref{thm: KY14: homotopy simple} is optimal when $d=2$.
\par
(ii)
For the natural map and stabilization in degree $2$, the target
second homotopy groups are both isomorphic to $\Z$ with trivial
fundamental-group action by the target decomposition in
Theorem \ref{thm: KY14: (m,n)=(3,1)} and Lemma \ref{lmm: the case d=3}.
The maps are surjective by Theorem \ref{thm: KY14: (m,n)=(3,1)}.
Equivariance therefore identifies each induced map on $\pi_2$, up
to sign, with the augmentation $\Lambda\to\Z$, whose kernel is
$(t-1)\Lambda$.
\par
(iii)
In degree $3$, Lemma \ref{lmm: the case d=3} and the target
decomposition in Theorem \ref{thm: KY14: (m,n)=(3,1)} give
\[
\begin{split}
\pi_3(\po^{3,3}_1(\R))&\cong\pi_3(S^2\vee S^3)\cong\Z\oplus\Z,\\
\pi_3((\Omega^2_3\CP^2)^{\Z_2})
&\cong\pi_5(S^5)\oplus\pi_4(S^2)\cong\Z\oplus\Z/2.
\end{split}
\]
The natural map is surjective on $\pi_3$ by the stability theorem
and cannot be injective. The stabilization in degree $3$ is also noninjective on $\pi_3$:
the natural map in degree $4$ induces an isomorphism on $\pi_3$,
so the target of stabilization has third homotopy group
$\Z\oplus\Z/2$. The stabilization is surjective on $\pi_3$
by the stability theorem, and thus has nonzero kernel.
Consequently neither part of Theorem \ref{thm: KY14: (m,n)=(3,1)}
can be strengthened to a homotopy equivalence through dimension $d$
for $d=2$ or $3$.
\qed
}
\end{remark}
\begin{remark}
{\rm 
Here is a second proof of Lemma~\ref{lmm: optimality d=2}, using
two closed subsets of the quadratic model. Translation of the variable
removes the linear term from the first polynomial and gives a homeomorphism
$\pomod{2}\cong\R\times\mathcal P_2$, where $\mathcal P_2$ consists
of triples of the form
\[
 f=z^2+c,\qquad (g,h)=zU+V,\qquad c\in\R,\quad U,V\in\R^2.
\]
Thus $\mathcal P_2$ is identified with the open set $Y$ of allowed
coefficients $(c,U,V)$. For $c<0$ both vectors
$V+\sqrt{-c}\,U$ and $V-\sqrt{-c}\,U$ must be nonzero.
For $c=0$ the condition is $V\ne0$. For $c>0$ the two roots
of $f$ are nonreal, and the condition is $(U,V)\ne(0,0)$.

Let $Y_-$ and $Y_+$ be the closed subsets of $Y$ given by $c\leq0$
and $c\geq0$, respectively, and put $Y_0=Y_-\cap Y_+$.
Both halves are manifolds with boundary $Y_0$. Their boundary
inclusions are cofibrations, and each interior is homotopy equivalent
to its half. On $c<0$ the two evaluations $V\pm\sqrt{-c}\,U$ give
independent coordinates for fixed $c$, so $Y_-\simeq S^1\times S^1$.
On $c>0$ the pair $(U,V)$ ranges over $\R^4\setminus\{0\}$, so
$Y_+\simeq S^3$. At $c=0$ the direction of $V$ gives $Y_0\simeq S^1$.

The two evaluations extend to $c=0$, where both equal $V$.
Under these equivalences the inclusion $Y_0\to Y_-$ therefore
represents the diagonal $\Delta:S^1\to S^1\times S^1$. The inclusion
$Y_0\to Y_+$ is null-homotopic since $\pi_1(S^3)=0$. Taking the
homotopy pushout gives $Y\simeq {\rm Cofib}(\Delta)\vee S^3$.
The torus automorphism $(z,w)\mapsto(z,wz^{-1})$ carries the diagonal
to $S^1\times\{1\}$. Collapsing this circle gives $S^1\vee S^2$ up to
homotopy, as in the proof of Lemma~\ref{lmm: optimality d=2}.
This proves \eqref{eq: po231}.
\qed
}
\end{remark}

\section*{Acknowledgements}
The authors would like to thank Professor Martin Guest for his many
valuable insights and suggestions, especially concerning Segal's classical work.

\par
The authors also acknowledge substantial mathematical assistance
from ChatGPT (OpenAI) in the preparation and revision of this paper.
In particular, ChatGPT identified a gap in an earlier argument for
homotopy stability and proposed the degree-filtration and
universal-cover homology arguments used to repair it (Theorem
\ref{thm: simple up to dimension} and Lemma
\ref{lmm: nilpotent covering homology}). It also contributed to the
revised proof of Lemma \ref{lmm: Lemma 3.3 in KY13}, suggested the use
of the Maslov index in Remark \ref{remark: space of equiv. maps}, and
proposed the geometric proof of Lemma \ref{lmm: optimality d=2}.
We also acknowledge the assistance of ChatGPT Pro in developing the proof
of the Arnold--Vassiliev case announced in the introduction.
The authors take full responsibility for the mathematical statements
and proofs.
\par
The second author was supported by 
JSPS KAKENHI Grant Number JP22K03283. 
This work was also supported by the Research Institute for Mathematical Sciences, a Joint Usage/Research Center located in Kyoto University.

 

\begin{thebibliography}{99}

\bibitem{Ar}
V. I. Arnold,
Certain topological invariants of algebraic functions, (Russian),
Trudy Moskov. Obshch., {\bf 21} (1970), 27-46
\bibitem{BHMM}
C. P. Boyer, J. C. Hurtubise, B. M. Mann and R. J. Milgram,
The topology of instanton moduli spaces, I: The Atiyah-Jones conjecture,
Ann. Math., {\bf 137}, (1993), 561-609.
\bibitem{BHMM2}
C. P. Boyer, J. C. Hurtubise, B. M. Mann and R. J. Milgram,
The topology of the space of rational maps into generalized flag manifolds, Acta Math., {\bf 173}, (1994), 61-101.
\bibitem{Br}R. W. Brockett,
Some geometric questions in the theory of linear systems,
IEEE Trans. Automat. Control, {\bf 21} (1976), 449--455.
\bibitem{Brown}K. S. Brown,
{\it Cohomology of Groups},
Graduate Texts in Mathematics, {\bf 87}, Springer-Verlag, New York, 1982.
\bibitem{CCMM}F. R. Cohen, R. L. Cohen, B. M. Mann and R. J. Milgram,
The topology of rational functions and spaces of divisors,
Acta Math., {\bf 166}, (1991), 163--221.
\bibitem{CMM}F. R. Cohen, M. Mahowald and R.J. Milgram,
The stable decomposition for the double loop space of a sphere,
Proc. Symp. Pure Math., {\bf 32} (1978), 225-228.
\bibitem{FW}B. Farb and J. Wolfson,
Topology and arithmetic of resultants, I,
New York J. Math., {\bf 22} (2016), 801-821.
\bibitem{GKY2}
M. A. Guest, A. Kozlowski and K. Yamaguchi,
Spaces of polynomials with roots of bounded
multiplicity,
Fund. Math., {\bf 116} (1999), 93--117.
\bibitem{Hatcher}A. Hatcher,
{\it Algebraic Topology}, Cambridge University Press, Cambridge, 2002.
\bibitem{James}I. M. James, Reduced product spaces,
Ann. Math.,  {\bf 62} (1955), 170-197.
\bibitem{JS}I. M. James and G. Segal,
On equivariant homotopy type, Topology, {\bf 17} (1978), 267-272.
\bibitem{KS}M. Konstantinov and J. Smith,
Monotone Lagrangians in $\CP^n$ of minimal Maslov number $n+1$,
Math. Proc. Cambridge Philos. Soc. {\bf 171} (2021), 1--21.
\bibitem{KY8}A. Kozlowski and K. Yamaguchi, 
The homotopy type of spaces of resultants of bounded multiplicity, Topology
Appl., {\bf 232} (2017), 112-139.
\bibitem{KY9}A. Kozlowski and K. Yamaguchi, 
The homotopy type of spaces of rational curves on a toric variety,
Topology Appl., {\bf 249} (2018), 19-42.
\bibitem{KY12}A. Kozlowski and K. Yamaguchi,
Spaces of non-resultant systems of bounded multiplicity determined by a toric 
variety, Topology Appl., {\bf 337} (2023), Paper no. 108626, 30p.
%
\bibitem{KY15}A. Kozlowski and K. Yamaguchi,
Spaces of non-resultant systems of real bounded multiplicity determined by a toric variety, 
J. Math. Soc. Japan, {\bf 77} (2025), 1137-1181.
\bibitem{KY13}A. Kozlowski and K. Yamaguchi,
Spaces of non-resultant systems of bounded multiplicity with real coefficients,
Algebraic \& Geometric Topology, {\bf 26} (2026), 2589--2633.
\bibitem{May2} J. P. May, Infinite loop space theory,
Bull. Amer. Math. Soc., {\bf 83} (1977), 456-494.
\bibitem{Milgram}R. J. Milgram, Iterated loop spaces, Ann. Math., {\bf 54} (1966), 386-403.
\bibitem{Se0}G. B. Segal,
Configuration spaces and iterated loop spaces, Invent. Math., {\bf 21}
(1973), 213-221.
\bibitem{Se}G. B. Segal,
The topology of spaces of rational functions,
Acta Math., {\bf 143} (1979), 39--72.
\bibitem{Va}V. A. Vassiliev,
Complements of discriminants of smooth maps, 
Topology and Applications, 
Amer. Math. Soc.,
Translations of Math. Monographs \textbf{98},
1992 (revised edition 1994).



\end{thebibliography}
\end{document}